\documentclass[11pt,reqno]{amsart}

\usepackage{amsmath,amsthm,amssymb,enumitem,xspace, calc}

\newcommand{\Z}{\mathbb{Z}}

\newcommand{\R}{\mathbb{R}}

\DeclareMathOperator{\conv}{conv}

\DeclareMathOperator*{\argmin}{arg\,min}

\DeclareMathOperator*{\proj}{proj}

\newcommand{\norm}[2]{\left\|#2\right\|_{#1}}

\newcommand{\partialDerivative}[2]{\frac{\partial#1}{\partial #2}}

\newcommand{\bigO}[1]{\ensuremath{\mathcal{O}(#1 )}}

\newtheorem{proposition}{Proposition}
\newtheorem{theorem}{Theorem}
\newtheorem{corollary}{Corollary}
\newtheorem{lemma}{Lemma}
\theoremstyle{definition}
\newtheorem{definition}{Definition}
\newtheorem{example}{Example}
\theoremstyle{remark}
\newtheorem{remark}{Remark}

\newcommand\defeq{\stackrel{\mathclap{\text{\tiny def}}}{=}}

\newcommand*{\ones}{\ensuremath{\mathbf{e}}}

\DeclareMathOperator{\meet}{\wedge}
\DeclareMathOperator{\join}{\vee}

\usepackage{algorithm,algorithmicx}
\usepackage{algpseudocode}
\usepackage{bbm}
\usepackage[mathlines,pagewise,left]{lineno}
\usepackage{amssymb}
\usepackage{color,xcolor,tikz}
\usetikzlibrary{arrows,shapes,chains}  
\usepackage{amsmath}
\usepackage{bcol}
\usepackage{multirow, array}
\usepackage{caption}
\usepackage{xspace}
\usepackage[numbers]{natbib}
\usepackage{comment}
\usepackage[toc,page]{appendix}
\usepackage{graphicx} 
\usepackage{booktabs,stackengine} 
\usepackage{pdflscape}
\usepackage[export]{adjustbox}
\usepackage[font=scriptsize]{subfig}
\usepackage{mathtools}
\usepackage{hyperref}
\usepackage{bm}
\usepackage[normalem]{ulem}

\usepackage{amsthm}
\theoremstyle{definition}

\theoremstyle{remark}
\usepackage[framemethod=TikZ]{mdframed}

\DeclareMathOperator*{\one}{\textbf{1}}

\newcommand{\rev}[1]{{\color{red} #1}}
\newcommand{\set}[1]{\ensuremath{\mathcal{#1}}}
\DeclareMathOperator*{\minimize}{\ensuremath{\mathbf{minimize}}}

\DeclareMathOperator*{\minimum}{\ensuremath{\mathbf{minimum}}}

\newcommand{\pivot}[3][k+1]{\ensuremath{\mathtt{pivot}^{#1}(#2,#3,\bm\alpha,\bm\beta)}}

\newcolumntype{\resetRow}{>{\global\let\currentrowstyle\relax}}
\newcolumntype{^}{>{\currentrowstyle}}

\def\SingleSpacedXI{\linespread{1.1}}
\SingleSpacedXI

\newcommand{\EO}{\ensuremath{\mathsf{EO}}\xspace}
\newcommand{\EG}{\ensuremath{\mathsf{EG}}\xspace}

\newcommand{\slow}{\ensuremath{\mathsf{Gurobi}}\xspace}

\newcommand{\fastp}{\ensuremath{\mathsf{Pivot}}\xspace}
\newcommand{\mip}{\ensuremath{\mathsf{MIP}}\xspace}
\newcommand{\sm}{\ensuremath{\mathsf{SFM}}\xspace}
\newcommand{\sparsity}{\ensuremath{\mathsf{Sparsity}}\xspace}
\newcommand{\bigM}{\ensuremath{\textsf{Big-M}}\xspace}
\newcommand{\strongMip}{\ensuremath{\textsf{Strong-MIP}}\xspace}

\usepackage[colorinlistoftodos,prependcaption,textsize=small]{todonotes}

\title[Convex Submodular Minimization with Indicator Variables]{Convex Submodular Minimization with Indicator Variables}
\author{Andr\'es G\'omez and Shaoning Han}
\thanks{ \noindent \hskip -5mm	
	A. G\'{o}mez: Daniel J. Epstein Department of Industrial and Systems Engineering, University of Southern California, CA 90089.   \texttt{gomezand@usc.edu}. \\
	S. Han: Department of Mathematics, and Institute of Operations Research and Analytics, National University of Singapore. Singapore 119076. \texttt{shaoninghan@nus.edu.sg}.
}
\begin{document}
\maketitle
\begin{center}
	June 2025
\end{center}
\begin{abstract}
	\vskip 3mm
	\noindent  We study a general class of convex submodular optimization problems with indicator variables. Many applications such as the problem of inferring Markov random fields (MRFs) with a sparsity or robustness prior can be naturally modeled in this form. We show that these problems can be reduced to binary submodular minimization problems, possibly after a suitable reformulation, and thus are strongly polynomially solvable. 
	 Furthermore, we develop a parametric approach for computing the associated extreme bases under certain smoothness conditions. This leads to a fast solution method, whose efficiency is demonstrated through numerical experiments.
	
	\vspace{1em}
\end{abstract}
\section{Introduction}
Given $\bm a,\bm d\in\R^n$, we consider the problem of the form
\begin{equation}\label{eq:intro}
	\minimize_{\bm{x}\in \R^n,\bm{z}\in\{0,1\}^n}\left\{f(\bm{x})-\bm a^\top \bm x+\bm{d^\top z}: \;\ell_iz_i\le x_i\le u_iz_i\,\forall i=1,\dots,n\right\}
\end{equation}
where:
\begin{enumerate}
	\item function $f:\R^{n}\to\R$ is convex and (continuous) submodular. 
	\item bounds $\bm{\ell}\in \underline{\R}^n$ and $\bm{u}\in \overline{\R}^n$ are possibly infinite, where $\underline\R\defeq\R\cup\{-\infty\}$ and $\overline{\R}\defeq\R\cup\{\infty\}$, and satisfy $\bm\ell\le \bm u$.
\end{enumerate}
Here we adopt the convention that $0\cdot(\pm\infty)=0$. Under this convention, if $z_i=0$, the constraints enforce $x_i=0$; if $z_i=1$,  then $x_i$ is activated and allowed to take any value in $[\ell_i,u_i]$, incurring a fixed cost $c_i$. Observe that we do not assume that $\ell_i\leq 0\leq u_i$ for any $i\in[n]$, and thus \eqref{eq:intro} is general enough to include the constraint where a continuous variable is either zero or bounded away from zero. The convex submodular term $f(x)$ can be used to capture the pairwise similarity or data fidelity of statistical models \cite{bertsimas2018predictive, bertsimas2024slowly}, which makes formulation \eqref{eq:intro} a natural choice for regression problems involving smoothness and combinatorial priors, such as sparse signal denoising and outlier detection in dynamic systems; see Section~\ref{sec:application} for a detailed discussion of applications. 

Submodular functions of binary variables are often equivalently represented as set functions characterizing the diminishing return property. They arise pervasively in combinatorial optimization \cite{nemhauser1978analysis,edmonds1970submodular}, with classical examples including cut capacity functions of networks \cite{megiddo1974optimal} and rank functions of matroids \cite{bergmann1929axiomatik,whitney1935abstract}, and are often associated with discrete optimization problems that admit efficient algorithms and theoretical guarantees \cite{fujishige2005submodular,iwata2008submodular}.  Recently, there has been growing interest in submodular optimization over continuous domains \cite{bach2019submodular, bach2012structured, hassani2017gradient, bian2017continuous, staib2017robust}, partially stimulated by applications in machine learning. However, less effort has been devoted to investigating submodular optimization problems involving both continuous and discrete variables \cite{yu2024constrained}. To the best of our knowledge, structured problems of the form~(1) have not been studied systematically in the literature.

Some special cases where $f(\bm x)$ is a quadratic function have been studied in literature. Observe that when $f(\bm x)=\bm x^\top \bm Q\bm x$, $\ell_i=-\infty$, $u_i=+\infty$ and $d_i=\lambda>0$ for all $i$, substituting out binary variables $z$ yields an equivalent unconstrained optimization problem
\[ \bm x^\top \bm Q\bm x-\bm a^\top \bm x+\lambda\norm{0}{\bm x}, \] where $\norm{0}{\bm x}$ denotes the so-called nonconvex $\ell_0$-``norm" and is defined as the number of nonzero components in $\bm x$. Without additional structure imposed over $Q$, this problem is in general $\mathcal{NP}$-hard as it subsumes the sparse linear regression problem as a special case. Notably, certain tractable cases emerge when $f(\bm x)=\bm x^\top \bm Q\bm x$ is convex and submodular, which is equivalent to $\bm Q$ being a \emph{Stieltjes} matrix, that is, $\bm Q$ is positive definite and $Q_{ij}\le 0$ for all $i\neq j$. In particular,  \citet{atamturk2018strong} show that if $\bm a\ge0$, $\ell_i=-\infty$ and $u_i=\infty$ for all $i$, then \eqref{eq:intro} can be recast as a binary submodular minimization problem, rendering it strongly polynomial solvable in theory. In addition, \cite{kim2003exact} indicates that under more restrictive conditions, \eqref{eq:intro} can be addressed via semidefinite programming. However, it remains an open question whether such polynomial solvability results can be extended to more general settings, allowing positive $a_i$ and finite $\ell_i$ and $u_i$.

Another important special case arises when $f(x)=\sqrt{\sigma^2+\sum_{i=1}^n c_ix_i^2}$, where $\sigma\ge0$ and $\bm c>\bm 0$. This form occurs widely in risk averse optimization, including mean risk minimization \cite{atamturk2008polymatroids}, Value-at-Risk minimization \cite{ghaoui2003worst}, and distributionally robust optimization \cite{zhang2018ambiguous}. Problem~\eqref{eq:intro} with the diagonal conic quadratic $f$ is first studied by \citet{atamturk2019lifted}. More recently, \citet{gomez2021strong} shows that when $u_i=\infty\;\forall i\in[n]$, the problem \eqref{eq:intro} admits an exact conic quadratic relaxation and is therefore polynomially solvable in these settings.

In practice, \eqref{eq:intro} can be solved using mixed-integer optimization (MIO) approaches.  On one hand, the natural relaxation obtained by relaxing $\bm z\in\{0,1\}^n$ to $[0,1]^n$ provides a lower bound of \eqref{eq:intro}. On the other hand, by fixing $\bm z$ to a specific binary vector, \eqref{eq:intro} reduces to a tractable convex optimization problem whose optimal value leads to an upper bound on the original problem \eqref{eq:intro}. Therefore, these bounds can be incorporated into black-box branch-and-bound algorithms for solving \eqref{eq:intro} exactly. In certain cases where $f(\bm x)=\bm x^\top \bm Q\bm x$ is a Stieltjes quadratic form, \citet{atamturk2021sparse} propose stronger conic relaxations by convexifying low-dimensional quadratic terms, which outperform the standard big-M relaxation. Similar ideas are also explored in solving general quadratic optimization with indicator variables \cite{gomez2021outlier, HGA2023, frangioni2020decompositions, shafiee2024constrained}. However, despite the potential advantages of MIO methods, they can suffer from scalability issues as the problem size grows. Our numerical experiments also confirm this point, highlighting the limitations of pure MIO approaches in large-scale settings.

\subsection*{Contributions} The contributions of this paper are two-fold. 

\subsubsection*{{1}. We show that if, for any fixed binary $\bm z\in\{0,1\}^n$, the corresponding box-constrained convex optimization problem derived from~\eqref{eq:intro} can be solved in (strongly) polynomial time, then the original mixed-integer submodular minimization problem~\eqref{eq:intro} is also (strongly) polynomially solvable. }\ \\

The result is established by introducing additional artificial binary variables and reducing \eqref{eq:intro} to minimizing a certain binary submodular function $v(\cdot)$  --a class of problems which admits polynomial time algorithms \cite{iwata2001combinatorial,iwata2009simple}, where each evaluation of $v(\cdot)$ relies on solving a box-constrained convex optimization involving $f(\cdot)$. In particular, when $f(\bm x)=\bm x^\top \bm Q\bm x$ is a Stieltjes quadratic form, our result implies that the corresponding mixed-integer quadratic optimization problem is strongly polynomially solvable, regardless of the sign of coefficients $\bm a$, thereby addressing the gap discussed above in the literature. Moreover, we further extend the results to non-Stieltjes quadratic objectives by leveraging the combinatorial structure of the matrix $\bm Q$.

\subsubsection*{{2}. We develop a fast method for computing extreme bases of the binary submodular function $v(\cdot)$ in question.}\ \\

Roughly speaking, an extreme base (the formal definition is given later in Definition~\ref{def:extremeBasis}) consists of $n+1$ evaluations of $v(\cdot)$, which are required in each iteration of all existing generic binary submodular minimization (BSM) algorithms. In our setting, computing these key quantities boils down to solving $n+1$ convex optimization problems, which can be expensive and renders solving \eqref{eq:intro} via BSM more conceptual than practical.  To overcome this bottleneck, we propose a parametric algorithm that  computes the extreme base progressively with a total computational cost comparable to a single evaluation of $v(\cdot)$. The proposed method offers benefits both theoretically and practically. First, it reduces the overall complexity of solving \eqref{eq:intro} by a factor $\bigO{n}$. Second, and more importantly, it makes solving \eqref{eq:intro} as a BSM problem practically feasible.  Experimental results show that our new method for solving \eqref{eq:intro} achieves an order-of-magnitude speedup over state-of-the-art MIO approaches, while also delivering superior solution quality.

\subsection*{Outline}  In \S\ref{sec:preliminaries} we introduce notations and necessary preliminaries for the paper. In \S\ref{sec:application}, we discuss applications of the mixed-integer optimization problem \eqref{eq:intro} in detail. In \S\ref{sec:polytime} we prove that \eqref{eq:intro} can be reduced to a binary submodular minimization problem and can be solved in polynomial time. We also discuss the extension of the result in quadratic cases. In \S\ref{sec:acceleration} we develop the parametric algorithm for computing the extreme bases of binary submodular functions and specialize it to quadratic and conic quadratic cases. In \S\ref{sec:computations}, we test the solution efficacy of the method proposed in this work on combinatorial Markov random field inference problems and present computational results. Finally, in \S\ref{sec:conclusion} we conclude the paper. 

\section{Preliminaries}\label{sec:preliminaries}
In this section we first introduce the concepts related to submodularity and notations used throughout the paper, and then briefly review the solution methods for binary submodular minimization (BSM) in literature. 

\subsection{Submodularity: definitions and notations.} Given an integer $n\in \Z_{++}$, we let $[n]\defeq\{1,\dots,n\}$.  We use bold symbols to denote vectors and matrices. For any $\bm x\in\R^n$, $\bm Q\in\R^{n\times n}$ and index sets $\alpha,\beta\subseteq[n]$, we denote by $\bm x_\alpha$ the subvector of $\bm x$ corresponding to the indices in $\alpha$, and $\bm Q_{\alpha\beta}$ the submatrix of $\bm Q$ with rows indexed by $\alpha$ and columns indexed by $\beta$.
We denote the vector of all zeros by $\bm{0}$ and the vector of ones by $\one$ (whose dimensions can be inferred from the context). 
 Given $i\in [n]$, we also let $\bm{e}^i$ be the $i$-th coordinate vector of $\R^n$.  
 We denote $\underline{\R}\defeq\R\cup\{-\infty\}$ and $\overline{\R}\defeq\R\cup\{\infty\}$ and we adopt the convention that $0\cdot(\pm\infty)=0$. For example, given decision variables $z\in \{0,1\}$ and $x\in \R$, constraint $-uz\le x\le uz$ with $u=\infty$ is equivalent to the complementarity constraint $x(1-z)=0$.  For a differentiable function $g:\R^n\to\R$ and $\alpha\subseteq [n]$, define $\nabla_\alpha g(\cdot, \bm x_{\alpha^c}):\R^{\alpha}\to\R^{\alpha}$ by $(\nabla_\alpha g(\bm x))_i=\frac{\partial}{\partial x_i}g(\bm x)\;\forall i\in\alpha$, where $\alpha^c$ is the complement of $\alpha$ in $[n]$. Additionally, if $g(\cdot)$ is strongly convex, then $\nabla_\alpha g(\cdot,\bm x_{\alpha^c})$ is invertible for any fixed $\bm x_{\alpha^c}$, and its inverse is denoted by $\nabla_\alpha^{-1}g(\cdot;\bm x_{\alpha^c})$.

Given two vectors $\bm{y}^1$ and $\bm{y}^2\in \R^n$, define the \emph{meet}
$\bm{y}^1\meet \bm{y}^2\in \R^n$ and the \emph{join} $\bm{y}^1\join \bm{y}^2\in \R^n$ to be the component-wise minimum and maximum of $\bm{y}^1$ and $\bm{y}^2$, respectively; we also define $\bm{y}^1\circ\bm{y}^2\in\R^n$ as the Hadamard (entrywise) product. By above notations, a set $\set L\subseteq \R^n$ is called a \emph{lattice} if any $\bm y^1,\bm y^2\in\set L$ implies that $\bm y^1\join \bm y^2$ and $\bm y^1\meet \bm y^2$ belong to $\set L$. A function $f:\R^n\to\R$ is \emph{submodular} over a lattice $\set L$ if for any $\bm y^1$ and $\bm y^2\in\set L$, one has $f(\bm y^1)+f(\bm y^2)\ge f(\bm y^1\meet \bm y^2)+f(\bm y^1\join \bm y^2).$ 
Proposition~\ref{prop:equiv-def-submodular} below provides several equivalent definitions of submodular functions; see \cite{topkis1978minimizing} for their reference.
\begin{proposition}[\citet{topkis1978minimizing}]\label{prop:equiv-def-submodular}
	The following statements hold true.\begin{itemize}
		\item (Zeroth-order definition) A function $f:\R^{n}\to\R$ is submodular if and only if for all $c_i,c_j>0,i\neq j$ and $\bm{y}\in\R^{n}$, it holds that \[f(\bm{y}+c_i\bm{e^i})+f(\bm{y}+c_j\bm{e^j})\ge f(\bm{y})+f(\bm{y}+c_i\bm{e^i}+c_j\bm{e^j}).\]
		
		\item (First-order definition) If $f:\R^{n}\to\R$ is differentiable, then $f$ is submodular if and only if $\partialDerivative{f}{y_i}(\bm{y}+c_1\bm{e^j})\le\partialDerivative{f}{ y_i}(\bm{y}+c_2\bm{e^j})$ for all $\bm y \in\R^n$, $i\neq j$ and $c_1\ge c_2$.
		\item (Second-order definition) If $f:\R^{n}\to\R$ is twice differentiable, then $f$ is submodular if and only if $\partialDerivative{^2 f(\bm{y})}{y_i\partial y_j}\le 0$ for all $\bm y\in\R^n$ and $i\neq j$.
	\end{itemize}
\end{proposition}
We list two special classes of submodular functions that are closely related to the applications considered in Section~\ref{sec:application}. From the second-order definition we find that any function of the form $f(\bm{x})=\bm{x^\top Q x}$ is (strongly) convex submodular if $\bm{Q}$ is a Stieltjes matrix, that is, $\bm{Q}\in\R^{n\times n}$ is symmetric positive definite and $Q_{ij}\leq 0$ for all $i\neq j$.  In addition, for any univariate convex function $g(\cdot)$, the function $h(x_1,x_2)\defeq g({x_1-x_2})$ is a composition of a convex function and a difference function, which can be easily verified to be convex and submodular over $\R^2$; see \cite{topkis1998}. We also point out that the sum of submodular functions is submodular, and the translation of a submodular function is submodular. 	

\subsection{Binary submodular minimization}
In this paper, we convert the optimization problem \eqref{eq:intro} to a binary submodular minimization  problem. Thus, we recall some necessary background on BSM. Given a binary submodular function $g:\set Z\to\R$, where $\set Z\subseteq\{0,1\}^n$ is a binary lattice, consider the binary submodular minimization problem $\min \{g(\bm z):\bm z\in\set Z\}$.

We treat a permutation over $[n]$ as a bijection $\pi:[n]\to[n]$, where $\pi_i\defeq \pi(i)\in[n]\;\forall i\in[n]$. Moreover, for any index set $\alpha\subseteq[n]$, we denote $\pi_{\alpha}\defeq\{ \pi_i:i\in\alpha \}$. In particular for $k\in\Z_{++}$, $\pi_{[k]}=\{\pi_1,\pi_2,\dots,\pi_k\}$. Denote the characteristic vector of index set $\alpha$ by $\ones^\alpha$, i.e. $\ones^{\alpha}_i=1$ if $i\in\alpha$ and $0$ otherwise. 
 Because there is a one-to-one correspondence between a binary vector $\bm z\in\{0,1\}^n$ and a subset $\alpha$ of $[n]$ through its characteristic vector, one often regards $g$ as a set function. We define $\Pi([n])$ as the set of permutations over $[n]$.  Extreme bases play an important role in BSM which we introduce as follows.
 \begin{definition}[Extreme base]\label{def:extremeBasis}
 	For any permutation $\pi\in\Pi([n])$, the extreme base $\bm y\in\R^{n}$ associated with $\pi$ is defined as 
 	\[ y_{\pi_i}=g\left( \ones^{\pi_{[i]}} \right)-g\left( \ones^{\pi_{[i-1]}} \right)\quad \text{for }i\in[n]. \]
 \end{definition}
 The computation of the extreme base induced by $\pi$ amounts to evaluating $\left\{ g\left(\ones^{\pi_{[i]}}\right) \right\}_{i=0}^{n}$. For convenience, we slightly abuse terminology and also refer to this sequence itself as the extreme base throughout the paper.
 
 Binary submodular minimization algorithms typically assume access to an \emph{evaluation oracle} for $g$. There are two main categories of approaches for BSM: combinatorial algorithms and convex optimization-based algorithms. The best combinatorial methods often enjoy a polynomial complexity in terms of an evaluation oracle \EO, where \EO denotes  the maximum amount of time required to evaluate $g(\ones^\alpha)$ for  $\alpha\subseteq[n]$. The seminal work of \citet{grotschel1981ellipsoid} introduced the first polynomial algorithm for BSM, with a strongly polynomial version later provided in \cite{grotschel1988geometric}. Other BSM combinatorial algorithms have also been developed subsequently in literature \cite{cunningham1985submodular,fleischer2000improved,fleischer2003push,schrijver2000combinatorial}. To the best of our knowledge, the current best complexity bound for general BSM is due to \citet{orlin2009faster}, whose algorithm runs in $\bigO{n^4\EG+n^7}$ time, where \EG stands for the maximum time of computing an extreme base. It is evident that $\EG$ is at most $n\cdot\EO$. In this paper, we will show that under mild conditions, one can achieve \EG=\EO for \eqref{eq:intro}. Although these combinatorial algorithms offer theoretical polynomial guarantees, they are often impractical due to high computational complexity. In fact, most of them have never been implemented.
 
 BSM can be converted to a convex optimization problem through \emph{Lov{\'a}sz extension}. More specifically, for any $\bm z\in\conv(\set Z)$, where $\conv(\set Z)$ is the convex hull of $\set Z$, the Lov{\'a}sz extension of $g$ at $\bm z$ is defined as the convex combination of the elements in the extreme base: $g^L(\bm z)\defeq\sum_{i\in[n]}(z_{\pi_i}-z_{\pi_{i+1}})g\left(\ones^{\pi_{[i]}}\right)+\left(1-z_{\nu_1}\right)g(\bm 0)$, where $\pi\in\Pi([n])$ is the permutation such that $z_{\pi_1}\ge\cdots\ge z_{\pi_n}$ and $z_{\pi_{n+1}}$ is defined as 0 for convenience. \citet{lovasz1983submodular} shows that $\min\{g(\bm z):\bm z\in\set Z  \}=\min\{g^L(\bm z):z\in\conv(\set Z)  \}$, where the latter problem is apparently convex. Convex optimization-based algorithms for BSM are more favorable than the combinatorial ones for practitioners, including cutting plane methods,  the minimum-norm point algorithm \cite{fujishige2006minimum, fujishige2011submodular}, and the conditional gradient method. We refer readers to Chapter~10 and Chapter~12 of the monograph \cite{bach2013learning} for a systematic treatment and experimental comparison of these approaches. Notably, all these methods require computing one extreme base in each iteration.  In the settings considered, evaluating $g$ is an expensive process as it requires solving a (convex submodular) minimization problem. We will develop a parametric algorithm to accelerate this process in Section~\ref{sec:acceleration}.

\section{Applications in MRF inference}\label{sec:application}
In this section, we begin by introducing Markov Random Field (MRF) inference problems and their applications across various domains. We then present two combinatorial variants of these problems and show how they can be reformulated in the form~\eqref{eq:intro}.

\begin{figure}[!h]
	\subfloat[1D]{\includegraphics[width=0.45\textwidth,trim={0cm 0cm 14cm 10cm},clip]{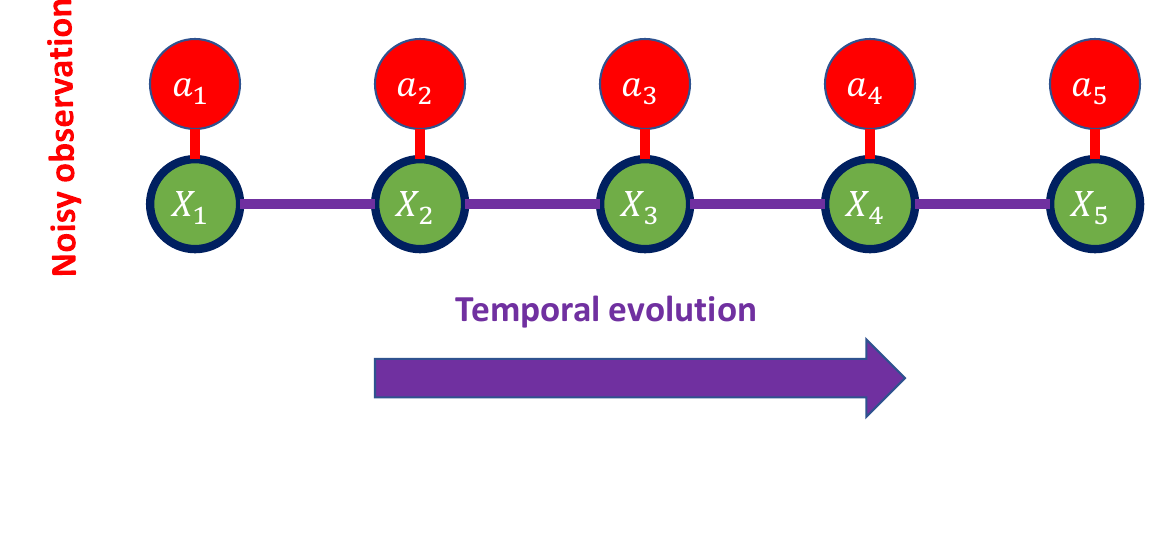}}\hfill\subfloat[2D]{\includegraphics[width=0.45\textwidth,trim={0cm 2cm 14cm 2cm},clip]{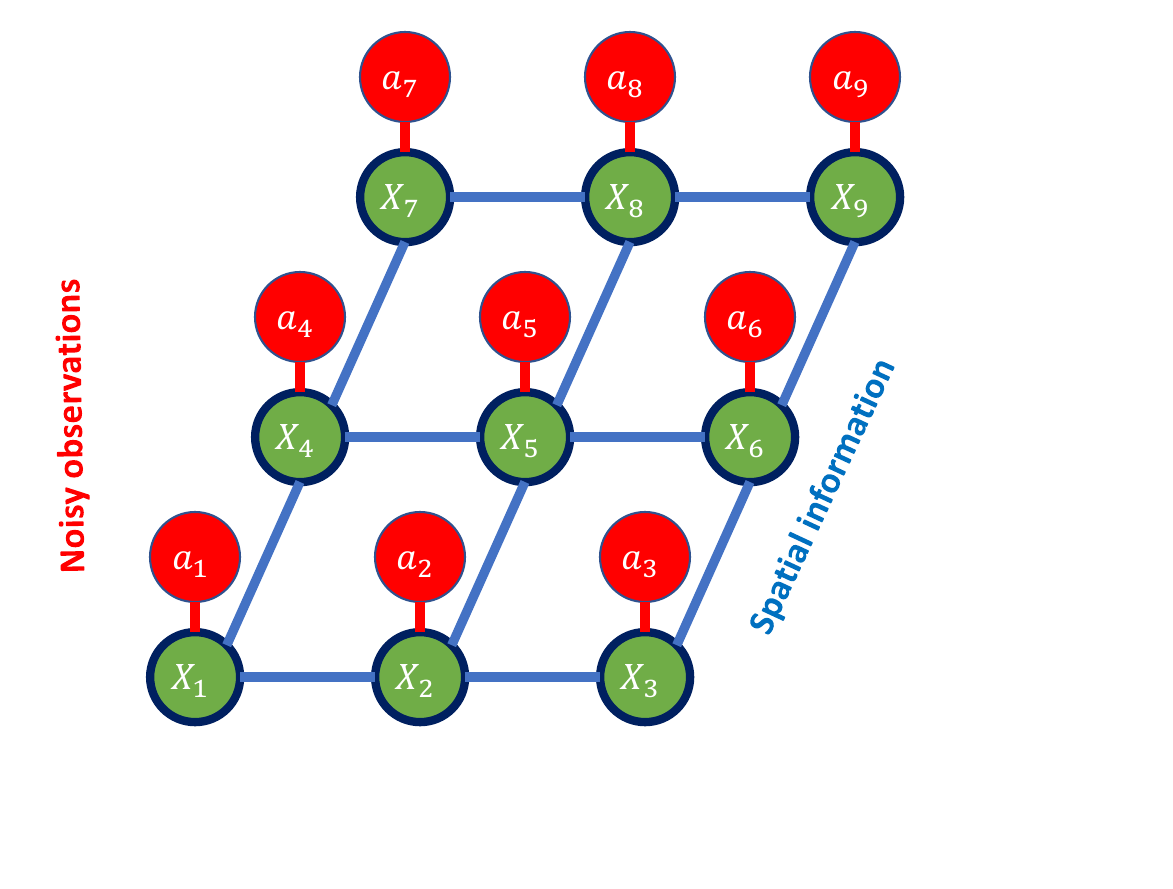}}\hfill	\newline
	\centering \subfloat[3D]{\includegraphics[width=0.75\textwidth,trim={0cm 4cm 2cm 0cm},clip]{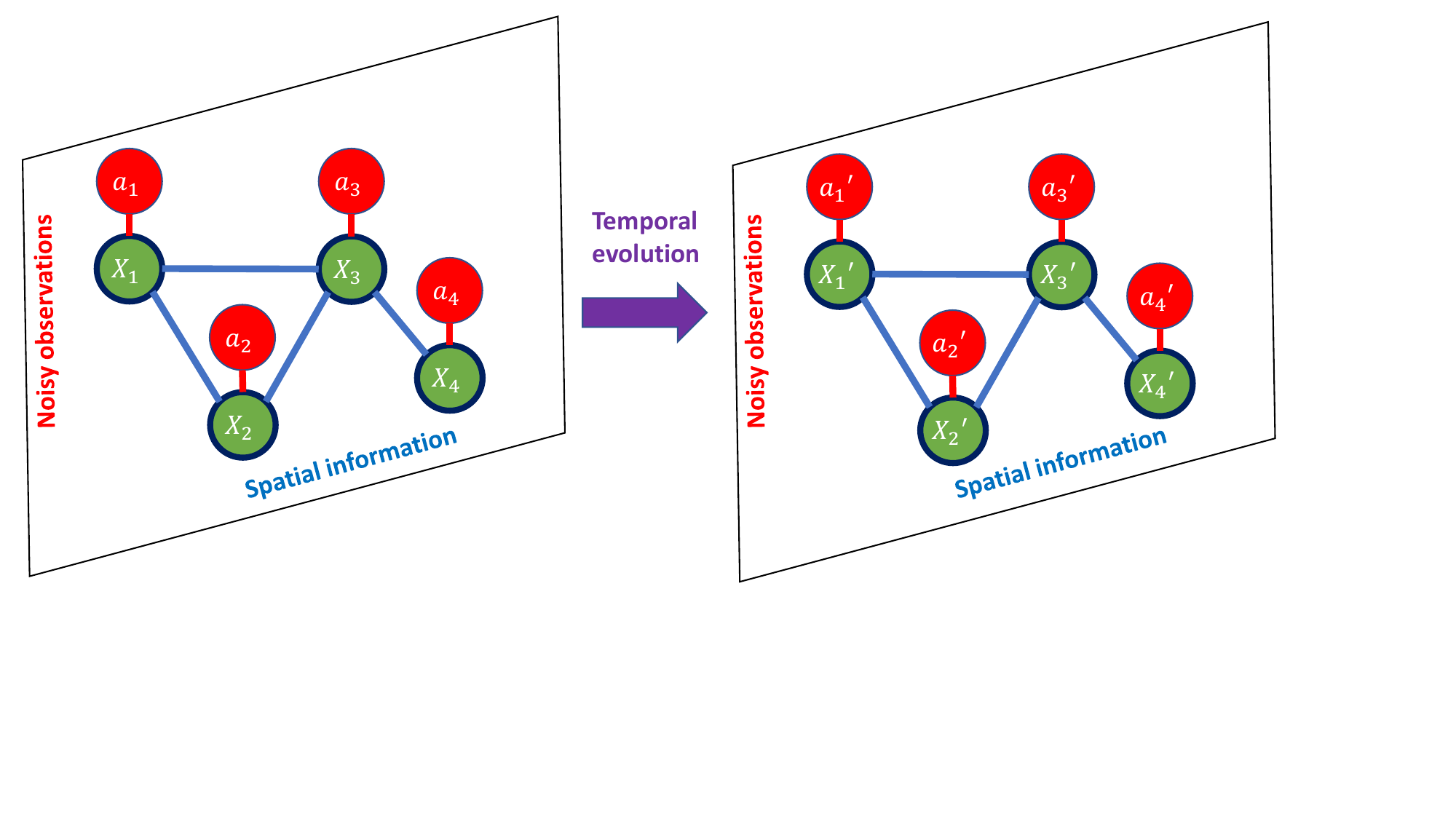}}
	\caption{Common topologies of MRFs, modeling spatial (blue) and temporal (purple) relationships. The true values of random variables $\bm{X}$ (green) are not observed directly and need to be inferred from the noisy observations $\bm{a}$ (red).}
	\label{fig:MRF}
\end{figure}

Markov random fields (MRFs) are popular graphical models pervasively used to represent spatio-temporal processes. They are defined on an undirected graph $\set G=(\set V,\set E)$, where there is random variable $X_i$ associated with each vertex $i\in \set V$. Each edge $[ i,j ] \in \set E$ represents the a relationship between the variables at their respective nodes $i$ and $j$; usually, these two variables should take similar values. Moreover, variables not connected by an edge are conditionally independent given realizations of all other variables. In the MRF inference problems we consider, noisy realizations $\{a_i\}_{i\in \set V}$ of the random variables $\bm{X}$ are observed, and the goal is to infer the true values of $\bm{X}$. Figure~\ref{fig:MRF} provides a depiction of this problem for three commonly-used structures of MRFs.

One-dimensional MRFs as depicted in Figure~\ref{fig:MRF} (A) are fundamental building blocks in time series analysis and signal processing \cite{angelov2006weighted,hochbaum2017faster,lu2022unified,mammen1997locally,restrepo1993locally,rinaldo2009properties}. They are typically used to model the evolution of a given process or signal over time. Two-dimensional MRFs as depicted in Figure~\ref{fig:MRF} (B) arise pervasively in image denoising \cite{boykov2006graph,boykov2001fast,hochbaum2001efficient,hochbaum2013multi,kolmogorov2004energy} and computer vision \cite{geman1986markov}. Each variable $X_i$ encodes the ``true" value of a pixel in an image, and edges encode the belief that adjacent pixels tend to have similar values. Two-dimensional MRFs also arise in ranking and selection problems based on similarity indexes~\cite{sun2019spectral,zhou2023sequential}. Three-dimensional MRFs as depicted in Figure~\ref{fig:MRF} (C) are used to model spatio-temporal processes \cite{fattahi2021scalable}. They are used in epidiomology \cite{besag1991bayesian,knorr1998modelling,morris2019bayesian} for example to track the spread of a disease over time.  In addition, MRFs over general graphs arise in semiconductor manufacturing \cite{ezzat2021graph,hochbaum2018adjacency}, bioinformatics \cite{eilers2005quantile}, criminology \cite{law2014bayesian}, spam detection \cite{hochbaum2019detecting}, among other applications.

\emph{Maximum a posteriori} estimates of the values of $\bm{X}$ can often be obtained as optimal solutions of the (continuous) MRF problem \cite{hochbaum2001efficient}
\begin{equation}\label{eq:MRF}
	\min\limits_{\bm{\ell}\leq \bm{x}\leq \bm{u}}\;\sum_{i\in \set V}h_i(x_i-a_i)+\sum_{[i,j]\in \set E}g_{ij}(x_i-x_j),
\end{equation}
where $h_i:\R\to\R_+$ and $g_{ij}:\R\to\R_+$ are appropriate convex nonnegative one-dimensional functions such that $h_i(0)=g_{ij}(0)=0$, and $\bm{\ell}\leq \bm{u}$ are (possibly infinite) lower and upper bounds, respectively, on the values of $\bm{X}$. From the comments following Proposition~\ref{prop:equiv-def-submodular}, it is clear that the objective of \eqref{eq:MRF} is a submodular function.  Functions $h_i$ and $g_{ij}$ are chosen depending on the prior distribution of the random variables and noise. Typically, functions $h_i$ are quadratic, corresponding to cases with Gaussian noise. The most common choices for functions $g_{ij}$ are absolute value functions $g_{ij}(x_i-x_j)=c_{ij}|x_i-x_j|$ with $c_{ij}\geq 0$, popular in statistics and signal processing \cite{sharpnack2012sparsistency,davies2001local} and referred to as total variation denoising problems, and quadratic functions $g_{ij}(x_i-x_j)=c_{ij}(x_i-x_j)^2$, in which case the graphical model is a Gaussian MRF (GMRF) and also corresponds to a Besag model \cite{besag1974spatial,besag1995conditional}. 

Clearly, problem \eqref{eq:MRF} is convex and can be solved using standard tools in the convex optimization literature. Specialized algorithms have also been proposed \cite{ahuja2004cut,hochbaum2001efficient}, whose complexity is strongly polynomial for the special cases of total variation and Besag models (see also \cite{hochbaum2013multi} and the references therein). In this paper, we study two combinatorial extensions of \eqref{eq:MRF}. The first extension corresponds to the situation where $\bm{X}$ is sparse or, more generally, is assumed to take a baseline value (e.g., corresponding to the background of an image or the absence of a disease) in most of its coordinates. In such cases, statistical theory calls for the imposition of an $\ell_0$ regularization to penalize variables that differ from the baseline value. The second extension corresponds to the situation where the noisy observations are corrupted by a few but potentially gross outliers. In such cases, statistical theory calls for the simultaneous removal of data identified as corrupted and solution of \eqref{eq:MRF}. Both extensions involve combinatorial decisions: which random variables differ from the baseline value, and which data points should be discarded. In some applications,  sparse and robust priors discussed above are incorporated  in the model simultaneously, e.g., \cite{yan2022real}. 

It is well known that linear regression, one of the simplest statistical estimation methods, becomes NP-hard with the inclusion of either sparsity \cite{natarajan1995sparse} or robustness \cite{bernholt2006robust} as described above. Thus, approaches in the literature resort to approximations of the combinatorial problems, heuristics, or expensive mixed-integer optimization approaches to solve the exact problems. \emph{In this paper we show that for the case of \eqref{eq:MRF}, the aforementioned combinatorial extensions can in fact be solved in polynomial time by a reduction to submodular minimization.} We point out that an immediate application of submodular minimization techniques \cite{orlin2009faster} results in runtime of $\mathcal{O}(n^5\cdot \EO)$, where $\EO$ is the complexity of solving problem \eqref{eq:MRF} -- resulting for example in strongly polynomial but impractical complexities of $\mathcal{O}(n^8)$ for the case of total variation and Besag models, but those runtime can likely be improved (we present such an improvement in this paper). Indeed, the discovery of a (strongly) polynomial time algorithm for a problem has typically been closely followed by highly efficient methods.


Next we formally define the two combinatorial extensions of problem \eqref{eq:MRF} discussed above --the sparse MRF inference problem and the robust MRF inference problem-- and their MIO formulations.

\subsection{ Sparse MRF inference} If the underlying statistical process $\bm{X}$ is known to be sparse (e.g., most pixels in an image adopt the background color, or the disease under study is absent from most locations), then a sparsity prior can be included in \eqref{eq:MRF}, resulting in problems of the form
\begin{subequations}\label{eq:SparseMRF}
	\begin{align}
		\min_{\bm{x}\in \R^{\set V},\;\bm{z}\in \{0,1\}^{\set V}}\;&\sum_{i\in \set V}h_i(x_i-a_i)+\sum_{[i,j]\in \set E}g_{ij}(x_i-x_j)+\sum_{i\in \set V}d_iz_i\\
		\text{s.t. }\;\;&\bm{\ell}\circ\bm{z}\leq \bm{x}\leq \bm{u}\circ\bm{z},
	\end{align}
\end{subequations}
where $\bm{d}\geq 0$ and binary variables $\bm{z}$ are used to indicate the support of $\bm{x}$ -- note that while solutions satisfying $z_i=1$ and $x_i=0$ are feasible, since $d_i\geq 0$ there always exists an optimal solution where $z_i=0$ if $x_i=0$. 
If all coefficients $d_i$ are equal, that is, $\bm{d}=\mu \one$ for some $\lambda\geq 0$, then in optimal solutions of \eqref{eq:SparseMRF} we have that $\sum_{i\in \set V}d_iz_i=\lambda\|\bm{x}\|_0$. 
Alternatively, if priors on the probabilities $p_i<0.5$ that variable $X_i$ is non-zero are available, then one can set $d_i\propto \ln((1-p_i)/p_i))$. Note that 
if $\bm{X}$ adopts a non-zero baseline value in most of its coordinates, the problem can be transformed into \eqref{eq:SparseMRF} through a change of variables.

Using MIO to model inference problems with sparsity is by now a standard approach in statistics and machine learning \cite{bertsimas2015or,bertsimas2016best,cozad2014learning,wilson2017alamo}. Most existing approaches focus on problems with quadratic functions -- probably due to the availability of powerful off-the-shelf MIO solvers capable of handling such functions. State-of-the-art methods revolve around the perspective relaxation \cite{akturk2009strong,frangioni2006perspective,gunluk2010perspective}: if $h_i(x_i-a_i)=(x_i-a_i)^2$, then we can replace such terms with the  reformulation $\hat h_i(x_i,z_i)=a_i^2-2a_ix_i+x_i^2/z_i$, where we adopt the following convention of division by $0$: $x^2/z=0$ if $x=z=0$, and $x^2/z=\infty$ if $x\neq 0$ and $z=0$. Indeed, this conic quadratic reformulation is exact if $z_i\in \{0,1\}$, but results in stronger continuous relaxations whenever $z_i$ is fractional. Tailored branch-and-bound algorithms \cite{hazimeh2021sparse}, approximation algorithms \cite{xie2020scalable} and presolving techniques \cite{atamturk2020safe} which exploit the perspective reformulation have been proposed in the literature. Finally, \citet{atamturk2021sparse} derive improved conic relaxations specific to problem \eqref{eq:SparseMRF} for the case of quadratic functions with $\bm{\ell}=\bm{0}$.  

Two special cases of \eqref{eq:SparseMRF} have been identified to be polynomial-time solvable. First, if graph $\set G$ is a path or a tree, then \eqref{eq:SparseMRF} can be solved via dynamic programming \cite{liu2022graph, bhathena2024parametric}. Second, all functions are quadratic, $u_i=\infty$ for all $i\in \set V$ and $\bm{a}\geq 0$, then \eqref{eq:SparseMRF} can be reformulated as a binary submodular problem \cite{atamturk2018strong} and thus be solved in polynomial time. 
\emph{In this paper, we show that such a submodular reformulation of \eqref{eq:SparseMRF} is always possible, regardless of the bounds, observations $\bm{a}$ or (convex) functions $h_i$ and $g_{ij}$.}

\subsection{Robust MRF inference}\label{sec:robust-mrf-inference} If the noisy observations $\bm{a}$ are corrupted by gross outliers, then the estimates resulting from \eqref{eq:MRF} can be poor. Classical robust estimation methods in statistics \cite{rousseeuw1984least,rousseeuw1987robust} call for the removal of outliers such that the objective \eqref{eq:MRF} is minimized, that is, solving the optimization problem
\begin{align}\label{eq:TrimmedMRF}
	\min\limits_{\bm{\ell}\leq \bm{x}\leq \bm{u},\;\bm{z}\in \{0,1\}^{\set V}}\;&\sum_{i\in \set V}h_i(x_i-a_i)(1-z_i)+\sum_{[i,j]\in \set E}g_{ij}(x_i-x_j)+\sum_{i\in \set V}d_iz_i,
\end{align}
where $z_i=1$ if and only if observation $i$ is discarded. Robust estimators such as \eqref{eq:TrimmedMRF} are, in general, hard to compute \cite{bernholt2006robust}.
In the context of least squares linear regression, the associated robust estimator is called the Least Trimmed Squares \cite{rousseeuw2006computing}, which is even hard to approximate \cite{mount2014least}. Exact optimization methods \cite{zioutas2005deleting,zioutas2009quadratic} rely on reformulations such as
\begin{subequations}\label{eq:TrimmedMRFReformulation}
	\begin{align}
		\min\limits_{\bm{x},\bm{z},\bm{w}}\;&\sum_{i\in \set V}h_i(x_i-w_i-a_i)+\sum_{[i,j]\in \set E}g_{ij}(x_i-x_j)+\sum_{i\in \set V}d_iz_i\\
		\text{s.t. }\;&w_i(1-z_i)=0\quad \forall i\in \set V\label{eq:TrimmedMRFReformulation_compl}\\
		&\bm{x}\in [\bm{\ell},\bm{u}]^{\set V},\;\bm{z}\in \{0,1\}^{\set V},\;\bm{w}\in \R^{\set V}.
	\end{align}
\end{subequations}
Indeed, since $h$ is nonnegative and $h(0)=0$, we find that if $z_i=1$, then $w_i=x_i-a_i$ in any optimal solution and the associated term vanishes; on the other hand, if $z_i=0$, then $w_i=0$ and $h_i(x_i-w_i-y_i)=h(x_i-a_i)$ as intended. Observe that problem \eqref{eq:intro} assumes each continuous variable is paired with an indicator.  This assumption is made without loss of generality. Indeed, in cases where some continuous variables do not have corresponding indicators like \eqref{eq:TrimmedMRFReformulation}, it is always possible to introduce an artificial binary variable $z_i$ with $d_i=0$ for each $w_i$ to transform the general problem into the form of \eqref{eq:intro}.  

Constraints \eqref{eq:TrimmedMRFReformulation_compl} are typically reformulated as big-M constraints; unfortunately, the ensuing continuous relaxation is trivial (e.g., $\bm{x}=\bm{0}$, $\bm{z}\to \bm{0}$, $\bm{w}=\bm{x}-\bm{a}$ in optimal solutions of the convex relaxations, and the objective value is almost $0$), thus the methods do not scale well. A stronger, big-M free, reformulation was proposed in \cite{gomez2021outlier} for the special case where $\set G$ is a path and all functions are convex quadratic.

Note that NP-hardness of robust estimators in general, and Trimmed Least Squares in particular, does not imply that \eqref{eq:TrimmedMRFReformulation} is NP-hard. \emph{In fact, we show in this paper that it is polynomial-time solvable for arbitrary convex functions $h_i$ and $g_{ij}$ and arbitrary graphs $\set G$.}

\section{Equivalence with binary submodular minimization}\label{sec:polytime}
In this section, we show that \eqref{eq:intro} can be reduced to a binary submodular minimization problem (under additional mild conditions). Our derivations are based on the fact that complementarity constraints preserve (to some degree) the lattice structure, and rely on the  following lemma.
\begin{lemma}[\citet{topkis1978minimizing}, Theorem 4.2]\label{lem:marginal-function}
	Given lattices $\set U$ and $\set W$, assume function $\phi:\set U\times\set W\to\R$ is submodular on a sublattice $\set L\subseteq \set U\times \set W$. If $$\psi(\bm{w})\defeq\min_{\bm{u}}\left\{ \phi(\bm{u},\bm{w}):(\bm{u},\bm{w})\in \set L\right\}>-\infty\quad\forall \bm{w}\in \set W,$$ then the marginal function $\psi$ is submodular on the lattice $\proj_{\bm{w}}(\set L)\defeq\{\bm{w}\in \set W:\exists \bm{u}\in \set U \text{ s.t. }(\bm{u},\bm{w})\in \set L\}$.
\end{lemma}

\subsection{General polynomiality results}

We first discuss the case where $\bm{\ell}\geq 0$, that is, $\bm{x}$ is nonnegative. Given $u\in \overline \R$, define \begin{equation}\label{eq:defLp}\set L_+=\{(x,z)\in\R\times\{0,1\}:\ell z\leq x\leq uz\}
 \end{equation}
\begin{lemma}\label{lem:lattice-plus}
	If $0\leq\ell\leq u$, then set $\set L_+$ is a lattice.
\end{lemma}
\begin{proof}
	Consider any $(x_1,z_1),(x_2,z_2)\in\set L_+$. It suffices to prove the case of $\min\{z_1,z_2\}$=0 since the other case where $\min\{z_1,z_2\}=1$ is trivial. If $\min\{z_1,z_2\}=0$, then $z_1=0$ or $z_2=0$, which implies $x_1=0$ or $x_2=0$. Since $0\leq x_1,x_2\le u$, one can deduce that $\min\{x_1,x_2\}=0$ and $\max\{x_1,x_2\}\leq u$; thus, $(x_1,z_1)\meet(x_2,z_2)\in\set L_+$ and $(x_1,z_1)\join(x_2,z_2)\in\set L_+$. Therefore, $\set L_+$ is a lattice.
\end{proof}

\begin{figure}[H]
	\subfloat[\parbox{0.8\linewidth}{Bounded  lattice $\mathcal{L}_+$}]{\includegraphics[width=0.45\textwidth, 
		clip]{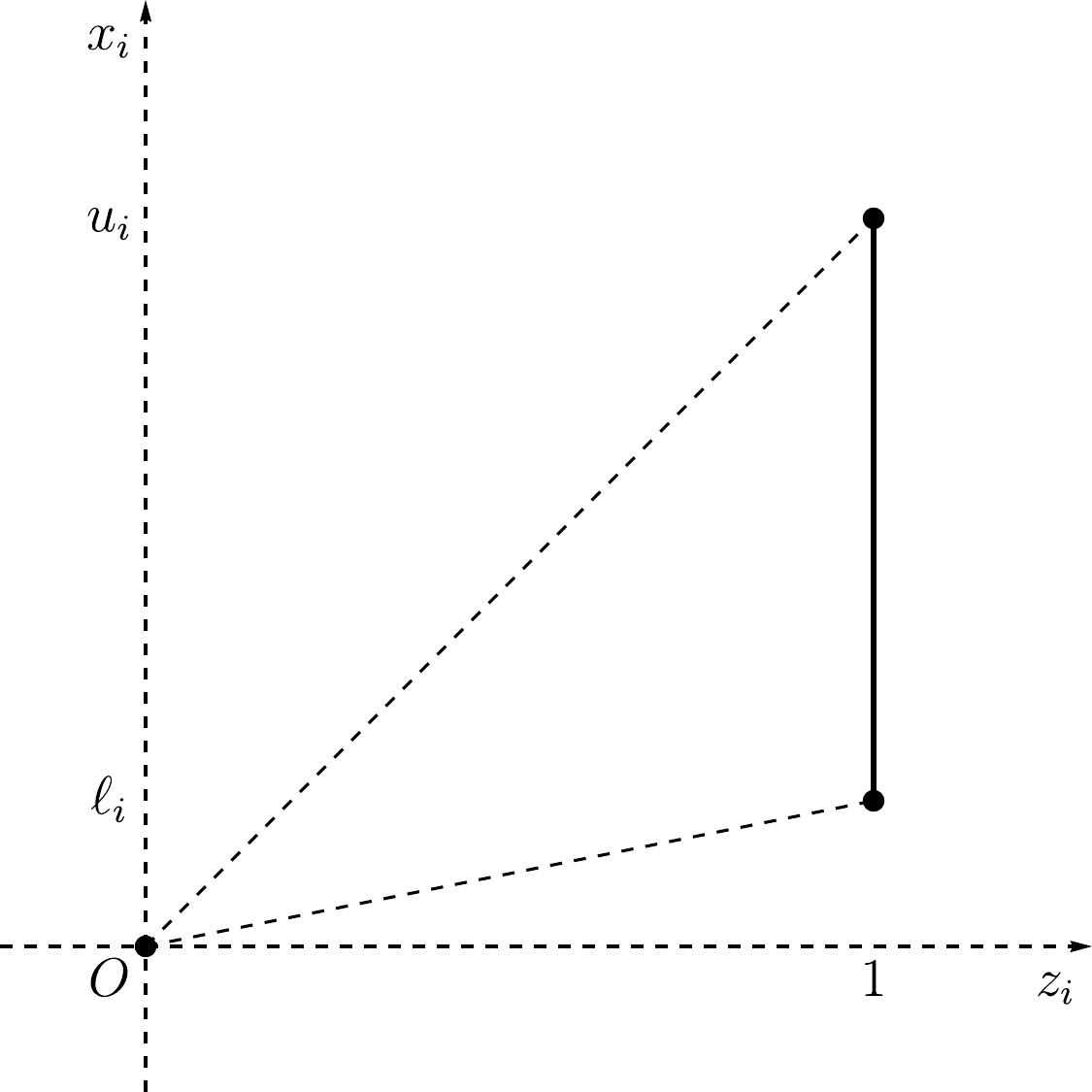}}\hfill\subfloat[{\makebox[0.8\linewidth][l]{Unbounded lattice $\mathcal{L}_+$}}]{\includegraphics[width=0.45\textwidth, clip]{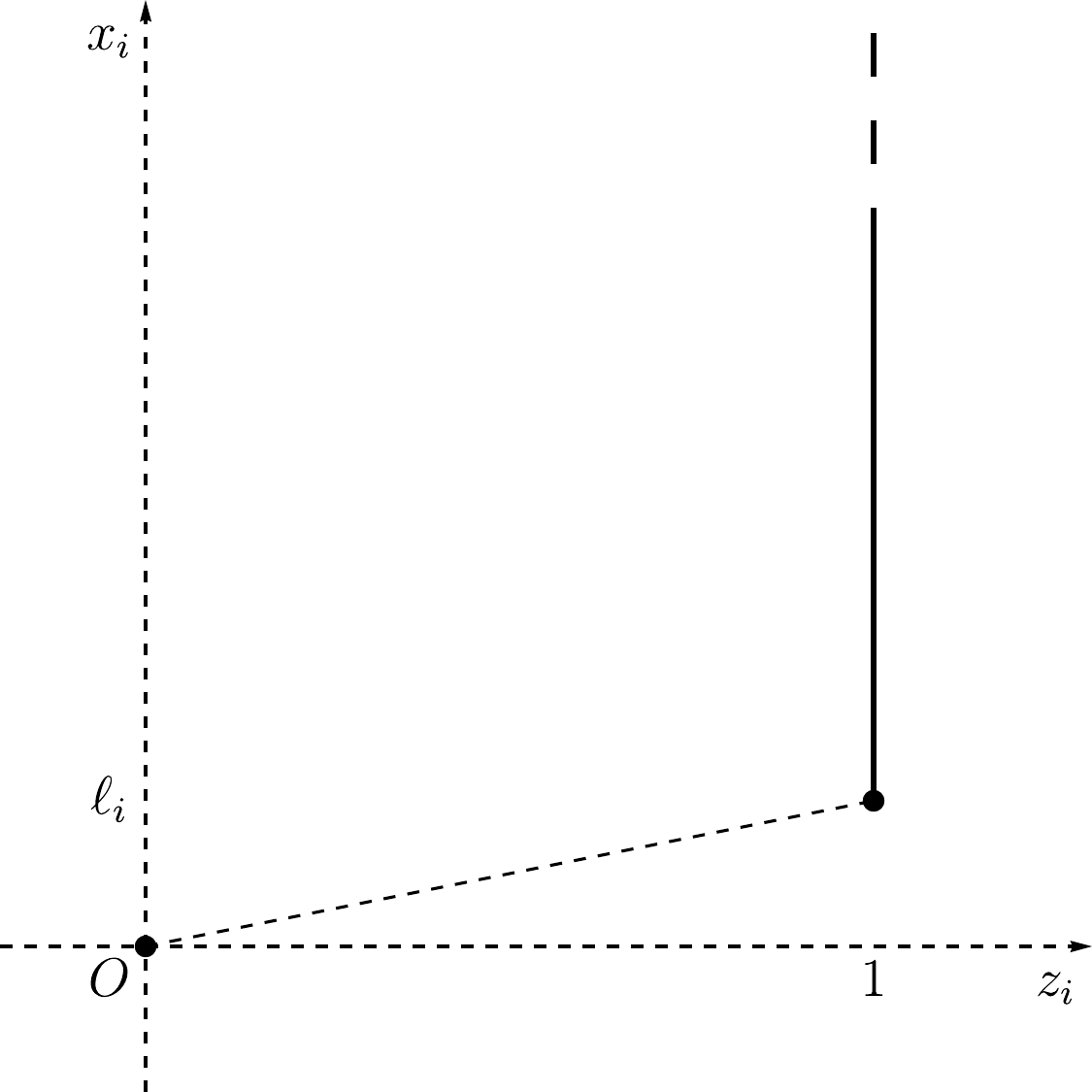}}\hfill\newline
	\subfloat[{\makebox[0.8\linewidth][l]{Not a lattice for $\ell_i<0<u_i$}}]{\includegraphics[width=0.45\textwidth, 
		clip]{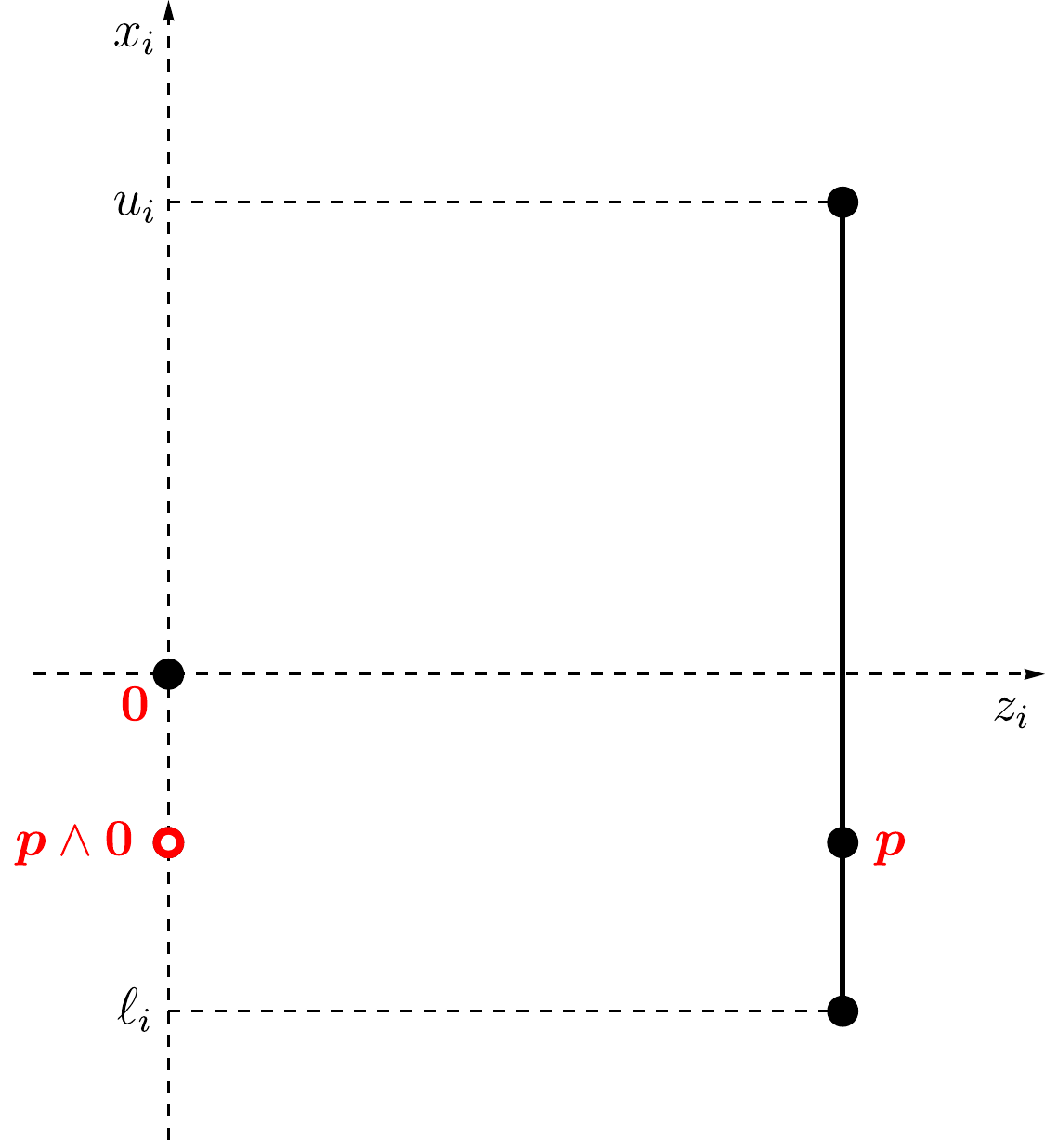}}\hfill\subfloat[{\makebox[0.8\linewidth][l]{Lifted feasible region $\mathcal{L}_\pm$}}]{\includegraphics[width=0.45\textwidth, trim = {2.8cm 4cm 3.2cm 0cm}, clip]{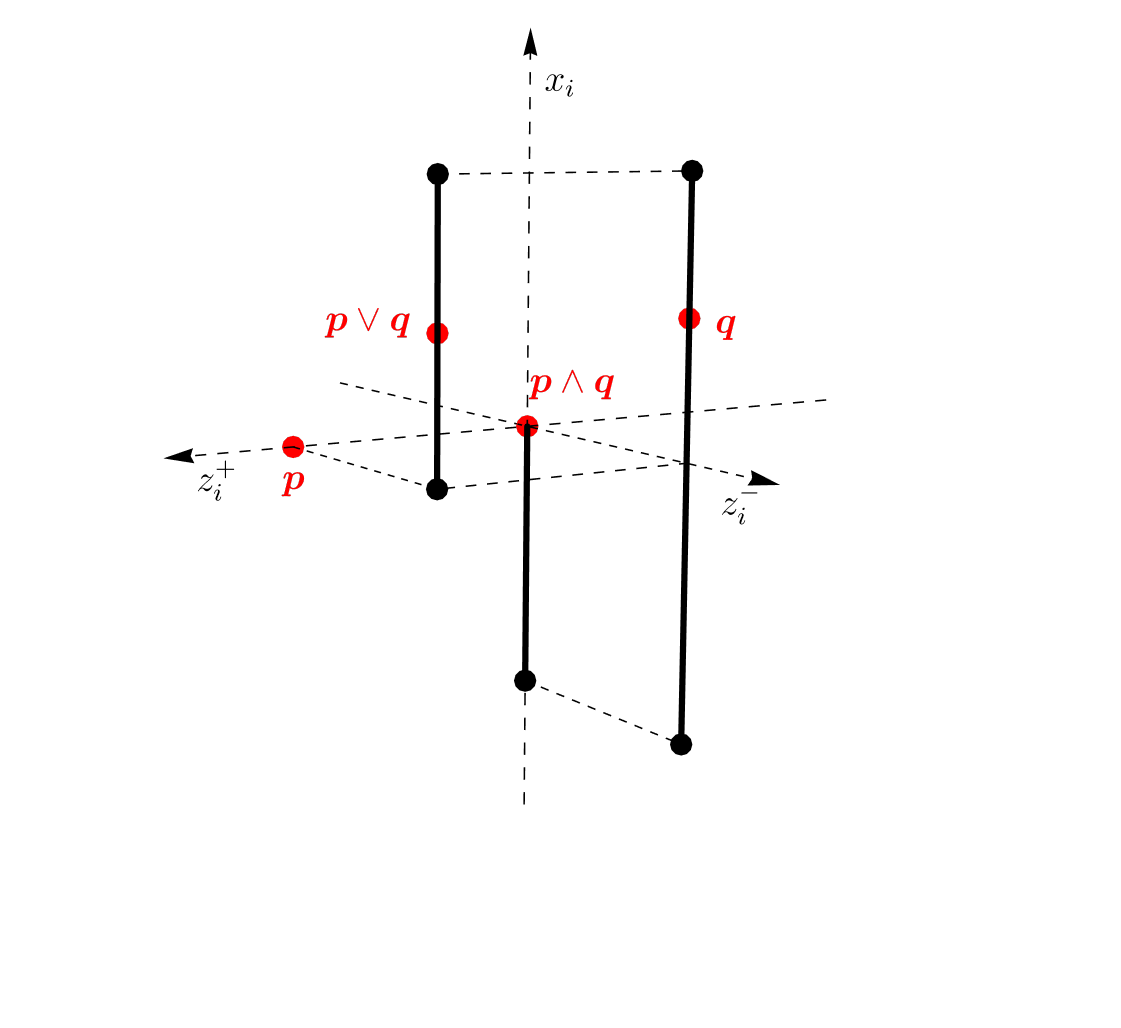}}
	\caption{Feasible region of mixed-integer submodular minimization problems}\label{fig:lattice}
\end{figure}

\begin{theorem}\label{thm:nonnegative-case}
	If $\bm{\ell}\in \R_+^n$, then the function \[ v_+(\bm{z})\defeq \min\;\left\{ f(\bm{x})-\bm a^\top \bm x:\,\bm{x}\in\R^n, \bm{\ell}\circ \bm{z}\leq \bm{x}\leq \bm{u}\circ{\bm{z}}\right\}\]
	is submodular on $\set Z$. 
\end{theorem}
\begin{proof}
	Note that the feasible region is a Cartesian product of $n$ lattices and thus is a lattice itself. The conclusion follows from Lemma~\ref{lem:marginal-function}. 
\end{proof}
If $\bm\ell\not\ge 0$, then the statement of Lemma~\ref{lem:lattice-plus} does not hold. Figure~\ref{fig:lattice}~(C) shows a counterexample where $\bm 0$ and $\bm p$ are feasible whereas their meet $\bm 0\meet\bm p$ is not. Consequently, function $v_+$ is not necessarily submodular. 
Next, we allow the continuous variables to be positive or negative and discuss how to address the non-lattice issue by expressing the feasible region in a lifting space.

As we show in Theorem~\ref{thm:general-case}, \eqref{eq:intro} can still be reformulated as a submodular minimization problem with the introduction of additional binary variables. Towards this goal, given $\ell\in \underline{\R}$ and $u\in \overline \R$, define additional sets 
\begin{align}\set\set \set L_-&=\{(x,z)\in\R\times\{0,1\}:\ell (1-z)\leq x\leq u(1-z)\} \label{eq:defLm}\\
\set L_\pm&=\{(x,z^+,z^-)\in\R\times\{0,1\}\times\{0,1\}:\ell (1-z^-)\leq x\leq uz^+\}.
	 \label{eq:defLpm}
\end{align}
\begin{lemma} \label{lem:lattice-gen}If $\ell\leq u\leq 0$, then $\set L_-$ is a lattice. If $\ell<0<u$, then $\set L_\pm$ is a lattice.
\end{lemma}

\begin{proof}
We prove just the result for $\set L_{\pm}$, as the proof of $\set L_-$ is analogous to the one of Lemma~\ref{lem:lattice-plus}. If $\ell\le0$ and $u\ge0$ are finite, then \[\set L_\pm=\left\{(x,z^+,z^-):-\ell z^-+\ell\le x\le uz^+ \right\}\cap\left(\R\times\{0,1\}^2\right)\] is a lattice as the intersection of two closed lattices is a closed lattice itself. In the general case where $\ell$ and $u$ are allowed to take infinite values, consider any $(x_1,z_1^+,z^-_1),(x_2,z_2^+,z_2^-)\in\set L_\pm$. Let $\hat u=x_1\join x_2\le u, \hat\ell=x_1\meet x_2\ge\ell$. Then $(x_i,z_i^+,z^-_i)\in\set{\hat L}_\pm:=\left\{(x,z^+,z^-)\in\R\times\{0,1\}^2:\hat\ell (1-z^-)\leq x\leq \hat uz^+\right\}$ for $i=1,2$. The conclusion follows from the lattice property of $\set{\hat L}_\pm$ and the inclusion $\set{\hat L}_\pm\subseteq\set L_\pm$.
\end{proof}

 To reformulate \eqref{eq:intro}, define $\set N_+\defeq\{i\in [n]:0\leq \ell_i\leq u_i\}$, $\set N_-\defeq\{i\in [n]:\ell_i\leq u_i\leq 0\}$ and $\set N_\pm\defeq\{i\in [n]:\ell_i<0< u_i\}$. For each $i\in \set N_{\pm}$ introduce binary variables $z_i^+=1$ if $x_i>0$ and $z_i^-=0$ if $x_i<0$, so that we can substitute $z_i=z_i^++(1-z_i^-)$ --note that we need to add constraint $z_i^-\geq z_i^+$ to rule out the impossible case where both $x_i>0$ and $x_i<0$. Figure~\ref{fig:lattice}~(D) shows the resulting lattice $\mathcal{L}_{\pm}$ from lifting the set in Figure~\ref{fig:lattice}~(C) through above transformation (without $z_i^-\ge z_i^+$). For convenience, for $i\in \set N_+$ we rename $z_i=z_i^+$ and for $i\in \set N_-$ we rename $z_i=1-z_i^-$. After performing the substitutions above, we find that \eqref{eq:intro} can be formulated as  
 \begin{subequations}\label{eq:submodular-reformulation}
 	\begin{align}
 		\min_{\bm{x},\bm{z}^+,\bm{z}^-}\;&f(\bm{x}) +\bm a^\top \bm x+\!\!\!\sum_{i\in \set N_+}\!\!d_iz_i^++\!\!\!\sum_{i\in \set N_-}\!\!d_i(1-z_i^-)+\!\!\!\sum_{i\in \set N_\pm}\!\!d_i(z_i^++1-z_i^-)\\
 		\text{s.t. }\;&\ell_i z_i^+\leq x_i\leq u_iz_i^+\quad\forall i\in \set N_+\label{eq:submodular-reformulation_const1}\\
 		&\ell_i (1-z_i^-)\leq x_i\leq u_i(1-z_i^-)\quad\forall i\in \set N_-\\
 		&\ell_i (1-z_i^-)\leq x_i\leq u_iz_i^+,\; z_i^-\geq z_i^+\quad\forall i\in \set N_\pm\\
 		&\bm{x}\in\R^n,\;\bm{z}^+\in\{0,1\}^{\set N_+\cup \set N_\pm},\;\bm{z}^-\in\{0,1\}^{\set N_-\cup \set N_\pm}.\label{eq:submodular-reformulation_const2}
 	\end{align}
 \end{subequations} 

\begin{proposition}\label{prop:latticeGen}
	The set defined by constraints
	\eqref{eq:submodular-reformulation_const1}-\eqref{eq:submodular-reformulation_const2} is a sublattice of $\R^n\times \{0,1\}^{\set N_+\cup \set N_\pm}\times \{0,1\}^{\set N_-\cup \set N_\pm}$.
\end{proposition}
\begin{proof}
	Each constraint involving $(x_i,z_i)$ jointly defines a lattice by Lemma~\ref{lem:lattice-plus} and Lemma~\ref{lem:lattice-gen}, and so does the constraint $z_i^- \geq z_i^+$.
\end{proof}

\begin{theorem}\label{thm:general-case}
	Function
	\begin{equation}\label{eq:subproblem_pm}
		  v_\pm(\bm{z}^+,\bm{z}^-)\defeq \min_{\bm{x}}\left\{f(\bm{x})-\bm a^\top \bm x:\eqref{eq:submodular-reformulation_const1}-\eqref{eq:submodular-reformulation_const2}\right\}
	\end{equation} 
	is submodular on $\{0,1\}^{\set N_+\cup \set N_\pm}\times \{0,1\}^{\set N_-\cup \set N_\pm}$. 
\end{theorem}
\begin{proof}
	Follows directly from Lemma~\ref{lem:marginal-function} and Proposition~\ref{prop:latticeGen}. 
\end{proof}
\begin{remark}
	In fact, Theorem~\ref{thm:nonnegative-case} and Theorem~\ref{thm:general-case} hold true for an arbitrary submodular (not necessarily convex) function $f$. However, if $f$ is not convex, the evaluation of the value function $v_{\pm}$ is in general not an easy task.
\end{remark}
\begin{remark}
	In some applications \cite{rosenbaum2010sparse, jia2015preconditioning, zhang2018sparse, gomez2019sign}, the artificial variables $z_i^+$ and $1-z_i^-$  are themselves important statistics to be inferred as they directly encode the sign of the covariates or signal in the underlying model. In such contexts, achieving \emph{sign consistency} (that is, correctly recovering the sign pattern of the true parameter vector) is often more important than mere support recovery.
\end{remark}
\begin{remark}\label{rem:delete}
	Observe that since $\bm{d}\geq 0$, constraints $z_i^-\geq z_i^+,\; \forall i\in \set N_\pm$ can be dropped from the formulation in principle. Indeed, if the constraints are removed and $z_i^+=1$, $z_i^-=0$ in an optimal solution of the resulting problem, then setting $z_i^+=0$ if $x_i\leq 0$ or $z_i^-=1$ if $x_i\geq 0$ results in a feasible solution with equal or better objective value. However, because $z_i^-\ge z_i^+$  shrinks the feasible region and reduces the number of cases to be considered in Section~\ref{sec:acceleration}, we retain them in the rest of the paper.  \qed
\end{remark}
\subsection{Implications for quadratic objectives}
	\citet{atamturk2018strong} show that problem  \[\min_{\bm{x}\in\R^n,\bm{z}\in \{0,1\}^n}\left\{\frac{1}{2}\bm x^\top \bm Q\bm x- \bm a^\top \bm x+\bm d^\top \bm z:\bm{x}\geq 0,\; x_i(1-z_i)=0\;\forall i\in [n]\right\}\]
	reduces to a submodular optimization problem provided that $\bm{Q}$ is a Stieltjes matrix and $\bm{a}\ge 0$. Theorem~\ref{thm:nonnegative-case} is a direct generalization, as it does not impose conditions on $\bm{a}$, allowing for arbitrary (nonnegative) variable lower bounds and arbitrary (finite or infinite) upper bounds on the continuous variables, and it holds for arbitrary (possibly non-quadratic) submodular functions. In this section, we generalize the polynomial solvability result to non-Stieltjes $\bm Q$, allowing more sign patterns of $\bm Q$ by exploiting graphical structures of the matrix.
	
	For a symmetric matrix $\bm Q\in\R^{n\times n}$, we denote by $G$ the (undirected) graph of $\bm Q$, where the set of vertices is $[n]=\{1,\dots,n\}$, and $i$ and $j$ are adjacent if and only if $i\neq j$ and $Q_{ij}\neq 0$. We denote by $G_-$ the graph on $\{1,2,\dots,n\}$ in which vertices $i$ and $j$ are adjacent if and only if $Q_{ij}<0$. The \emph{contraction} of an edge $e=(u,v)$ of a graph $G$ is to delete the edge $e$ and then identify its ends $u$ and $v$. Note that $G_-$ is a subgraph of $G$. We denote by $G/G_-$ the resulting graph by contracting all edges of $G_-$ in $G$.
\begin{theorem}\label{thm:miqp}
	Given a symmetric semidefinite matrix $Q\in\R^{n\times n}$ and the associated graphs $G$ and $G_-$, if $G/G_-$ is a bipartite graph, then the mixed-integer optimization problem 
	\begin{equation}\label{eq:miqp}
		\min\;\left\{\frac{1}{2}\bm x^\top \bm Q\bm x -\bm a^\top \bm x +\bm d^\top \bm z: \bm \ell\circ \bm z\le \bm x\le \bm u\circ \bm z,\;\bm z\in\{0,1\}^n \right\}.
	\end{equation}
is strongly polynomially solvable for all $\bm \ell\in\underline{\R}^n$ and $\bm u\in\overline{\R}^n$.
\end{theorem}
\begin{proof}
	We prove the result by properly changing signs of $x_i$ and reducing it to the Stieltjes case. Because $G/G_-$ is obtained by edge contraction, one can treat each vertex of $G/G_-$ as a subset of $[n]$. Moreover, since $G/G_-$ is bipartite, the vertices of $G/G_-$ can be partitioned into two parts $\mathcal{U}$  and $\mathcal{V}$ such that each edge of $G/G_-$ has one end in $\mathcal{U}$ and one end in $\mathcal{V}$. Define $\bar U =\cup_{U\in\mathcal{U}}U\subseteq[n]$, $\bar V=\cup_{V\in\mathcal{V}}V\subseteq[n]$, and a diagonal matrix $\bm D\in\R^{n\times n}$ with $D_{ii}=1$ if $i\in\bar U$ and $-1$ if $j\in\bar V$. Changing variables $\bm y=\bm D\bm x\Leftrightarrow 
\bm x=\bm{Dy}$, the problem \eqref{eq:miqp} is equivalent to
	\[ \min\;\left\{\frac{1}{2}\bm y^\top \bm D\bm Q\bm D\bm x -(\bm D\bm a)^\top \bm y +\bm d^\top \bm z: \bm \ell\circ \bm z\le \bm D\bm y\le \bm u\circ \bm z,\;\bm z\in\{0,1\}^n \right\}.  \]
	It remains to prove that $\bm{\bar Q}=\bm D\bm Q\bm D$ is a Stieltjes matrix, that is, $\bar Q_{ij}=D_{ii}D_{jj}Q_{ij}\le 0\,\forall i\neq j\in[n] $. Consider any $i,j\in[n]$ and $i\neq j$. If $Q_{ij}<0$, then $i$ and $j$ are identified in $G/G_-$ which implies either $i,j\in\bar U$ or $i,j\in\bar V$. In both cases, $D_i=D_j$, implying $\bar Q_{ij}<0$. If $Q_{ij}>0$, then one can deduce that either $i\in\bar U, j\in\bar V$ or $i\in\bar V$, $j\in\bar U$. In both cases, $D_iD_j=-1$, implying $\bar Q_{ij}<0$. This completes the proof.
\end{proof}
Note that  if $\bm Q$ is a Stieltjes matrix, then $G/G_-$ is a singleton. Therefore, Theorem~\ref{thm:miqp} includes the case of Stieltjes $\bm Q$ as a special case. Moreover, since the edge contraction of a tree always yields another tree, which remains bipartite, we can immediately deduce the following corollary.
\begin{corollary}
	If the graph $G$ of $\bm Q$ is a tree, then \eqref{eq:miqp} is strongly polynomially solvable.
\end{corollary}
We remark that when $G$ possesses specific structures, specialized algorithms may exist that are more efficient than solving \eqref{eq:miqp} via general submodular minimization.   In particular, when $G$ is a tree or even a path, \cite{liu2023graph} and \cite{bhathena2024parametric} show that \eqref{eq:miqp} can be solved in $\bigO{n^3}$  using dynamic programming approaches.
\section{Fast computation of extreme bases}\label{sec:acceleration}
In this section, we show that it is possible to compute an extreme base in the same complexity as a single evaluation of $v$ using a parametric algorithm, ultimately reducing the complexity of minimization algorithms by a factor of $n$.

In particular, we focus on a special case of \eqref{eq:intro}, for which we rename the indicator variables for the sake of algorithmic development, 
\begin{subequations}\label{eq:submodular-minimization-s}
	\begin{align}		
		\min_{\bm x\in\R^n,\, \bm s\in\{0,1\}^{2n-n_1}}\;&f(\bm x)-\bm a^\top \bm x+\bm d^\top \bm s\\
		\text{s.t. }\;	
		&\ell_is_i\le x_i\le u_is_i&i=1,\dots,n_1 \label{eq:submodular-s-positive-variable}\\
		&[x_{n_1+i}]_-s_{n_1+i}=0&i=1,\dots, n-n_1\\
		&[x_{n_1+i}]_+(1-s_{n+i})=0&i=1,\dots, n-n_1\label{eq:submodular-s-free-variable}\\
		&s_{n_1+i}\ge s_{n+i},&i=1,\dots, n-n_1\label{eq:constraint-s}
	\end{align}
\end{subequations}
where $f:\R^n\to\R$ is a convex submodular function, $\bm \ell\in\R_+^{n_1}$ and $\norm{\infty}{\bm u}<\infty$. The optimization problem \eqref{eq:submodular-minimization-s} corresponds to the case in \eqref{eq:submodular-reformulation} where $\mathcal{N}_+=[n_1]$, $\mathcal{N}_-=\emptyset$, and $\ell_i=-\infty$ and $u_i=\infty$ for all $i\in\mathcal N_{\pm}=\{n_1,\dots, n\}$.  Observe that we concatenate binary variables $\bm z_{\mathcal{N_+}}$, $\bm z_{\mathcal{N}_\pm}^-$ and $\bm z_{\mathcal{N}_\pm}^+$ in order into a single vector $\bm s$ here. The solution method to be proposed for solving \eqref{eq:submodular-minimization-s} in this section can be readily extended to more general cases including the case of $\mathcal{N}_-\neq\emptyset$ or bounded $\bm x_{\mathcal{N}_\pm}$-variables.

According to Theorem~\ref{thm:general-case}, problem \eqref{eq:submodular-minimization-s} can be polynomially solved by minimizing
\begin{equation}\label{eq:binary-submodular-minimization-s}
	\min\left\{ v(\bm s)+\bm d^\top \bm s:\; \bm s\in\{0,1\}^{2n-n_1},\,s_{n_1+i}\ge s_{n+i}\;\forall i=1,\dots, n-n_1 \right\},
\end{equation} 
where $v(\bm s)$ is a binary submodular function defined by
\begin{equation}\label{eq:subproblem-s}
	v(\bm s)=\min_{\bm x\in\R^n}\left\{ f(\bm x)-\bm a^\top \bm x:\; \eqref{eq:submodular-s-positive-variable}-\eqref{eq:submodular-s-free-variable}\right\}.
\end{equation}
For ease of exposition of the parametric approach to be proposed, throughout this section we additionally assume that $f$ is strongly convex and differentiable to ensure the finiteness of $v(\bm s)$ and uniqueness of the solution to subproblems \eqref{eq:subproblem-s}. However, we must point that these assumptions are not strictly required. If strong convexity fails to hold, one needs to first verify whether $v(\bm s)>-\infty$ for all feasible $\bm s$, which can be accomplished by verifying $v(\bm 1)>-\infty$. Secondly, multiple minima may exist for a given variable of $\bm s$, and the solution path needs to choose between these optimal solutions. If the differentiability of $f$ fails to hold, then one can use  subdifferential or directional derivative in place of $\nabla f$, and the analysis carried out in this section still holds up to moderate modification.

The workhorse behind all existing BSM algorithms is the efficient computation of extreme bases, which amounts to solving  a sequence of subproblems of the form \eqref{eq:subproblem-s} in our setting. The proposed method traces a solution path of minima of \eqref{eq:subproblem-s} as $\bm s$ varies in an isotonic manner.  More specifically, given a permutation $\pi:[2n-n_1]\to[2n-n_1]$ that is \emph{compatible} with \eqref{eq:constraint-s}, that is $s_{k}\ge s_{\ell}$ if $k=\pi_{n_1+i}$ and $\ell=\pi_{n+i}$ for any $i\in[n-n_1],$ our goal is to evaluate  functions $v(\ones^{\pi_{[i]}})$ progressively for all $i=0,\dots,2n-n_1$. Suppose that we have already computed $v(\ones^{\pi_{[k]}})$ for a certain $k=0,\dots,2n-n_1-1$ and let $\bm{\bar x}^{k}$ denote the associated optimal solution to \eqref{eq:subproblem-s}. We aim to evaluate $v(\ones^{\pi_{[k+1]}})$ and its optimal solution $\bm{\bar x}^{k+1}$, using $\bm{\bar x}^{k}$ as a warm start. Given $k\in\{0,1,\dots, 2n-n_1-1\}$, define function $v^\pi_{k+1}:\R\to\R$ as 
\begin{subequations}
	\begin{align*}
		v^\pi_{k+1}(y)=\min_{\bm x\in\R^{2n-n_1}}\;&f(\bm x)-\bm a^\top \bm x\\
		\text{s.t. }&x_{\pi_{k+1}}=y\\\tag{$P_{k+1}(y)$}\label{eq:subproblem-y}
		&\bm s=\ones^{\pi_{[k]}}\\
		&\ell_is_i\le x_i\le u_is_i&i\in[n_1]\backslash\{ \pi_{k+1} \}\\
		&[x_{n_1+i}]_-s_{n_1+i}=0&i\in[n-n_1]\backslash\{ \pi_{k+1}-n_1 \}\\
		&[x_{n_1+i}]_+(1-s_{n+i})=0&i\in[n-n_1]\backslash\{ \pi_{k+1}-n \}.
	\end{align*}
\end{subequations}
We denote by \eqref{eq:subproblem-y} the optimization subproblem defining $v_{k+1}^\pi(y)$. Moreover, denote the optimal solution to \eqref{eq:subproblem-y} by $\bm x^{{k+1}}(y)$. Here, we omit the dependence of $\bm{\bar x}^k,\;\bm x^{k}(y)$ and \eqref{eq:subproblem-y} on  permutation $\pi$.  The following lemma shows that $\bm x^{k+1}(\cdot)$ is isotonic in parameter $y$. 
\begin{lemma}\label{lem:isotonicity}
	The following statements hold true
	\begin{enumerate}[leftmargin=*]
		\item $\bm x^{{k+1}}({\bar x}^{k}_{\pi_{k+1}})=\bm{\bar x}^{k}$.
		\item $v(\ones^{\pi_{[k+1]}})=\min\limits_{y\in\mathcal{I}_{k+1}}v^\pi_{k+1}(y)$, where $\mathcal{I}_{k+1}=[\ell_{\pi_{k+1}},u_{\pi_{k+1}}]$ if $\pi_{k+1}\in[n_1]$, $\mathcal{I}_{k+1}=\{0\}$ if $\pi_{k+1}\in[n]\backslash[n_1]$, and $\mathcal{I}_{k+1}=[0,\infty)$ if $\pi_{k+1}\in[2n-n_1]\backslash[n]$.
		\item $\bm x^{k+1}(y_1)\ge \bm x^{k+1}(y_2)$ if $y_1\ge y_2$.
	\end{enumerate}
\end{lemma}
\begin{proof}
	The first two statements follow from the definition of the notations. Part~(3) is proved in Theorem~6.3, \cite{topkis1978minimizing}.
\end{proof}
Note that Lemma~\ref{lem:isotonicity} (1) and (2) hold for a generic function $f$. Submodularity is only used in part (3).  Lemma~\ref{lem:isotonicity}) brings the insights that we trace the path of all optimal solutions $\bm x^{k+1}(y)$ as $y$ is increased from ${\bar x}^k_{\pi_{k+1}}$ to reach the interval $\mathcal{I}_{k+1}$, and then as $y$ varies over $\mathcal{I}_{k+1}$. During this process, as implied by part (3) of Lemma~\ref{lem:isotonicity}, $\bm x^{k+1}(y)$ moves from $\bm{\bar x}^{k}$ to $\bm{\bar x}^{k+1}$ accordingly. To formally describe this routine, we introduce some index sets to represent the state of each coordinate $x^{k+1}_i(y)$. For a given $\bm x\in\R^{2n-n_1}$, define
\begin{equation}\label{eq:index-sets-x}
	\begin{aligned}
		&\alpha_0(\bm x)=\{ \pi_i\in[n_1]:\;\pi_i\notin\pi[k+1]\}\\
		&\underline{\alpha}(\bm x)=\{ \pi_i\in[n_1]:\;x_{\pi_i} =\ell_{\pi_i},\pi_i\in\pi[k]\}\\
		&\alpha_+(\bm x)=\{ \pi_i\in[n_1]:\;\ell_{\pi_i}< x_{\pi_i}  < u_{\pi_i},\pi_i\in\pi[k]\}\\
		&\overline{\alpha}(\bm x)=\{ \pi_i\in[n_1]:\;x_{\pi_i}  =u_{\pi_i},\pi_i\in\pi[k]\}\\
		&\beta_-(\bm x)=\{\pi_i\in[n]\backslash[n_1]:\;x_{\pi_i} <0,\pi_i\notin\pi[k+1]\}\\
		&\beta_\ominus(\bm x)=\{ \pi_i\in[n]\backslash[n_1]:\;x_{\pi_i} =0,\pi_i\notin\pi[k+1]\}\\
		&\beta_0(\bm x)=\{ \pi_i\in[n]\backslash[n_1]:\;\pi_i\in\pi[k],\pi_i+n-n_1\notin\pi[k+1] \}\\
		&\beta_\oplus(\bm x)=\{ \pi_i\in[n]\backslash[n_1]:\;x_{\pi_i} =0,\pi_i+n-n_1\in\pi[k] \}\\
		&\beta_+(\bm x)=\{ \pi_i\in[n]\backslash[n_1]:\;x_{\pi_i} >0,\pi_i+n-n_1\in\pi[k] \}\\
		&\gamma(\bm x)=\alpha_+(\bm x)\cup\beta_-(\bm x)\cup\beta_+(\bm x)\\
		&\gamma_0(\bm x)=\alpha_0(\bm x)\cup\beta_0(\bm x)\cup\beta_\ominus(\bm x)\cup\beta_\oplus(\bm x).
	\end{aligned}
\end{equation}
Intuitively, set $\alpha_0(\bm{x})$ is the set of variables in $[n_1]$ that have not yet been considered when computing $v_{k+1}^\pi$, and are fixed to $0$. Set $\overline{\alpha}(\bm x)$ are variables that reached their upper bound; since the path of solutions is isotonic, these variables will remain constant throughout the rest of the procedure. Sets $\underline{\alpha}(\bm x)$ and $\alpha_+(\bm x)$ are variables that may increase as the solution path is traced, and values of $y$ causing such variables to adopt a value different from the lower bound or reaching the upper bound for the first time need to be identified. For indices $i\in [n]\setminus[n_1]$, sets  $\beta_-(\bm{x})$ and $\beta_\ominus(\bm x)$ contain variables where both indicator variables (controlling lower and upper bounds) are set to zero, thus variables are non-positive; we distinguish between those that are strictly negative (and may increase as $y$ increases) and those that reached the value $0$. Set $\beta_0(\bm{x})$ contains variables where the indicator variable controlling the lower bound is set to $1$, and the indicator controlling the upper bound is $0$; those variables are simply fixed to $0$. Sets $\beta_\oplus(\bm x)$ and $\beta_+(\bm{x})$ are non-negative variables in the current iteration. Finally, set $\gamma(\bm{x})$ involves all variables not set to a bound (and thus increase continuously) if $y$ increases, and set $\gamma_0(\bm{x})$ contains all zero variables.

When it is clear from the context, the dependence on $\bm x$ of the index sets defined above will be omitted. Note that problem~\ref{eq:subproblem-y} is essentially a box-constrained convex optimization problem, allowing us to characterize its optimal solution $\bm x^{k+1}(y)$ in terms of KKT conditions. One equivalent variant of KKT conditions is stated in Lemma~\ref{lem:KKT} and no proof is needed. 
\begin{lemma}\label{lem:KKT}
	A point $\bm x\in\R^{2n-n_1}$ is  the optimal solution to \eqref{eq:subproblem-y} if and only if there exist index sets $\alpha_0$, $\underline{\alpha}$, $\alpha_+$, $\overline{\alpha}$, $\beta_-$, $\beta\ominus$, $\beta_0$, $\beta_\oplus$, $\beta_+$ and $\gamma=\alpha_+\cup\beta_-\cup\beta_+$ such that 
	\begin{equation}\label{eq:index->x}
		\begin{aligned}
			&x_{\pi_{k+1}}=y,\;\bm x_{\alpha_0}=\bm 0,\bm x_{\underline{\alpha}}=\bm \ell_{\underline{\alpha}},  \bm x_{\overline{\alpha}}=\bm u_{\overline{\alpha}},\\
			&\bm x_{\beta_\ominus}=\bm 0,\bm x_{\beta_\oplus}=\bm 0, \bm x_{\beta_0}=\bm 0, 
			\bm x_{\gamma}=\nabla^{-1}_\gamma f(\bm a_\gamma;\bm \ell_{\underline{\alpha}},\bm u_{\overline{\alpha}},\bm 0_{\gamma_0}, y)
		\end{aligned}		
	\end{equation}
	 satisfying
	 \begin{subequations}\label{eq:KKT-index}
	 	\begin{align}
	 		& \bm x_{\alpha_+}=\left(\nabla_\gamma^{-1} f(\bm a_{\gamma};\bm \ell_{\underline{\alpha}},\bm u_{\overline{\alpha}},\bm 0_{\gamma_0},y)\right)_{\alpha_+}\ge\bm \ell_{\alpha_+} \label{eq:KKT-alpha-lb}\\
	 		& \bm x_{\alpha_+}=\left(\nabla_\gamma^{-1} f(\bm a_{\gamma};\bm \ell_{\underline{\alpha}},\bm u_{\overline{\alpha}},\bm 0_{\gamma_0},y)\right)_{\alpha_+}\le\bm u_{\alpha_+}\label{eq:KKT-alpha-ub}\\
	 		&\bm x_{\beta_-}=\left(\nabla_\gamma^{-1} f(\bm a_{\gamma};\bm \ell_{\underline{\alpha}},\bm u_{\overline{\alpha}},\bm 0_{\gamma_0},y)\right)_{\beta_-}\le\bm 0\label{eq:KKT-beta-minus}\\
	 		&\bm x_{\beta_+}=\left(\nabla_\gamma^{-1} f(\bm a_{\gamma};\bm \ell_{\underline{\alpha}},\bm u_{\overline{\alpha}},\bm 0_{\gamma_0},y)\right)_{\beta_+}\ge\bm 0\label{eq:KKT-beta-plus}\\
	 		&\nabla_{\underline{\alpha}}f\left(\bm \ell_{\underline{\alpha}},\bm x_\gamma,\bm u_{\overline{\alpha}},\bm 0_{\gamma_0},y\right)\ge \bm a_{\underline{\alpha}}\label{eq:KKT-underline-alpha}\\
	 		&\nabla_{\overline{\alpha}}f\left(\bm \ell_{\underline{\alpha}},\bm x_\gamma,\bm u_{\overline{\alpha}},\bm 0_{\gamma_0},y\right)\le \bm a_{\overline{\alpha}}\label{eq:KKT-overline-alpha}\\
	 		&\nabla_{\beta_{\ominus}}f\left(\bm \ell_{\underline{\alpha}},\bm x_\gamma,\bm u_{\overline{\alpha}},\bm 0_{\gamma_0},y\right)\le \bm a_{\beta_{\ominus}}\label{eq:KKT-beta-ominus}\\
	 		&\nabla_{\beta_{\oplus}}f\left(\bm \ell_{\underline{\alpha}},\bm x_\gamma,\bm u_{\overline{\alpha}},\bm 0_{\gamma_0},y\right)\ge \bm a_{\beta_{\oplus}}\label{eq:KKT-beta-oplus}.	 		
	 	\end{align}
	 \end{subequations}
\end{lemma}
After substituting out $\bm x_\gamma$, one finds that variable sets $\alpha$'s and $\beta$'s are determined by inequalities \eqref{eq:KKT-index}, which are only related to parameter $y$ and problem data.  Hence,  Lemma~\ref{lem:KKT} reduces  the task of tracing the solution path $\bm x^{k+1}(y)$ to tracking the changes of $\alpha$'s and $\beta$'s. We now discuss the upper-bound condition of $y$ under which $\alpha$'s and $\beta$'s remain unchanged. As $y$ is increasing, $\bm x^{k+1}(y)$ is never decreasing by Lemma~\ref{lem:isotonicity}. This implies that constraints \eqref{eq:KKT-alpha-lb} and \eqref{eq:KKT-beta-plus} cannot block the increase of $y$.  On the other hand, because $f$ is a submodular function, one has $\frac{\partial^2}{\partial x_i\partial x_j}f(\bm x)\le0$, implying that the left-hand side of \eqref{eq:KKT-underline-alpha}--\eqref{eq:KKT-beta-oplus} is nonincreasing in $x_\gamma$ and $y$. Thus, \eqref{eq:KKT-overline-alpha} and \eqref{eq:KKT-beta-ominus} cannot block the increase of $y$ as well. Consequently, the next time that $\alpha$'s and $\beta$'s are altered can only happen when one of \eqref{eq:KKT-alpha-ub}, \eqref{eq:KKT-beta-minus}, \eqref{eq:KKT-underline-alpha} and \eqref{eq:KKT-beta-oplus} becomes active, depending on which of them occurs first. The corresponding threshold value of $y$ leads to a \emph{breakpoint} in the solution path and can be calculated by comparison. We formalize the above procedure and its conclusion as Algorithm~\ref{algo:pivot} and Proposition~\ref{prop:pivoting}, respectively.Define the tuple of index sets $\bm\alpha(\bm x)=(\alpha_0(\bm x),\underline{\alpha}(\bm x),\alpha_+(\bm x),\overline{\alpha}(\bm x))$ and $\bm\beta(\bm x)=(\beta_-(\bm x),\beta_{\ominus}(\bm x),\beta_0(\bm x),\beta_{\oplus}(\bm x),\beta_+(\bm x))$.
\begin{proposition}\label{prop:pivoting}
	Algorithm~\ref{algo:pivot} returns the next breakpoint $y_{\texttt{bp}}$ in the segment of the solution path starting from $y_0$ and ending with $\bar y$ if  one exists, and correctly updates the index sets $\bm\alpha(\bm x^{k+1}(y_{\texttt{bp}}))$ and $\bm\beta(\bm x^{k+1}(y_{\texttt{bp}}))$. Otherwise, it returns the end point $\bar y$ along with the original input index sets.
\end{proposition}

\begin{algorithm}[H]
	\caption{$\mathtt{pivot}^{k+1}(y_0,\bar y,\bm\alpha^0,\bm\beta^0)$ 
	}\label{algo:pivot}
	\footnotesize
	\begin{algorithmic}[1]
		\State $y_{\texttt{bp}}\gets y_0$, $\bm\alpha$$\gets$$\bm\alpha^0$, $\bm\beta$$\gets\bm\beta^0$, $\gamma\gets\alpha_+\cup\beta_-\cup\beta_+$, $\bm r\gets \bm{\infty}\in\overline{\R}^{2n-n_1}$ \Statex\Comment{$\bm{r}$ stores the potential thresholds of $y$}
		\For{$i\in\alpha_+$}
		\State $r_i\gets$ the root of $ \left(\nabla_\gamma^{-1} f(\bm a_{\gamma};\bm \ell_{\underline{\alpha}},\bm u_{\overline{\alpha}},\bm 0_{\gamma_0},y)\right)_{i}-u_i=0$
		\EndFor
		\For{$i\in\beta_-$}
		\State $r_i\gets$ the root of $ \left(\nabla_\gamma^{-1} f(\bm a_{\gamma};\bm \ell_{\underline{\alpha}},\bm u_{\overline{\alpha}},\bm 0_{\gamma_0},y)\right)_{i}=0$
		\EndFor
		\For{$i\in\underline{\alpha}\cup\beta_\oplus$}
		\State $r_i\gets$ the root of $\nabla_{i}f\left(\bm \ell_{\underline{\alpha}},\nabla_\gamma^{-1} f(\bm a_{\gamma};\bm \ell_{\underline{\alpha}},\bm u_{\overline{\alpha}},\bm 0_{\gamma_0},y),\bm u_{\overline{\alpha}},\bm 0_{\gamma_0},y\right)- a_{i}=0$ \label{line:lb-ri}
		\EndFor
		\State $i^*\gets\argmin\{r_i:i\in[2n-n_1]\}$ \Comment{Break the tie arbitrarily if any} \label{line:ratio test}
		\State $y_{\texttt{bp}}\gets\min\{ r_{i^*},\bar y \}$
		\If{$y_{\texttt{bp}}< \bar y$}
		\If{$i^*\in\alpha_+$}
		\State $\alpha_+\gets \alpha_+\backslash \{i^*\}$, $\overline{\alpha}\gets\overline{\alpha}\cup\{i^*\}$
		\EndIf
		\If{$i^*\in\beta_-$}
		\State $\beta_-\gets \beta_-\backslash \{i^*\}$, $\beta_\ominus\gets\beta_\ominus\cup\{i^*\}$
		\EndIf
		\If{$i^*\in\underline{\alpha}$}
		\State $\underline{\alpha}\gets \underline{\alpha}\backslash \{i^*\}$, $\alpha_+\gets\alpha_+\cup\{i^*\}$
		\EndIf
		\If{$i^*\in\beta_\oplus$}
		\State $\beta_\oplus\gets\beta_\oplus\backslash \{i^*\}$, $\beta_+\gets\beta_+\cup\{i^*\}$
		\EndIf
		\State $\gamma\gets\alpha_+\cup\beta_-\cup\beta_+$
		\EndIf
		\State \Return $($$y_{\texttt{bp}}$, $\bm\alpha$, $\bm\beta$, $\gamma$$)$
	\end{algorithmic}
\end{algorithm}
\begin{remark}
	Strictly speaking, the index sets $\alpha$'s and $\beta$'s generated during the execution of Algorithm~\ref{algo:pivot} do not exactly match those defined in \eqref{eq:index->x} and \eqref{eq:KKT-index}, because there may exist multiple index set configurations satisfying  \eqref{eq:index->x} and \eqref{eq:KKT-index} at breakpoints. Nonetheless, this does not affect the correctness of the algorithm.
\end{remark}
All root-finding equations arising in Algorithm~\ref{algo:pivot} are in terms of $y$. When any one of them has no solution, it is understood that the corresponding $r_i\gets\infty$. Denote the output of Algorithm~\ref{algo:pivot} by $\mathtt{pivot}^{k+1}(y_0,\bar y,\bm\alpha^0,\bm\beta^0)$. We now present Algorithm~\ref{algo:active-set}, which repeatedly calls Algorithm~\ref{algo:pivot} as a subroutine to compute the extreme bases $v(\ones^{\pi_{[k]}})$ for all values of $k$. 

\begin{algorithm}[H]
	\footnotesize
	\caption{Procedure to compute $\left\{v(\ones^{\pi_{[k]}})\right\}_{k=0}^{2n-n_1}$} 
	\label{algo:active-set}
	\begin{algorithmic}[H]
		\State  \textbf{Setup.} $\bm x\gets\argmin\{f(\bm x)-\bm a^\top \bm x:\,\bm x_{[n_1]}=\bm 0,\,\bm x_{[n]\backslash[n_1]}\le\bm 0\}$, $v(\bm 0)\gets f(\bm x)-\bm a^\top \bm x$, $\alpha_0\gets [n_1]$, $\beta_-\gets\{i: x_i<0  \}$,  $\beta_\ominus\gets[n]\backslash(\alpha_0\cup\beta_-)$,  $\underline{\alpha},\alpha_+,\overline{\alpha},\beta_0,\beta_\oplus,\beta_+\gets\emptyset$, $\gamma\gets\beta_-$ 
		\For{$k=0,1,\dots,2n-n_1-1$}
		\State  $y_0\gets x_{\pi_{k+1}}$ \Comment{Lemma~\ref{lem:isotonicity} Part~(1)} \label{line:start of the loop}
		\If{$\pi_{k+1}\le n_1$}
		\State $\alpha_0\gets \alpha_0\backslash[\pi_{k+1}]$
		\While{$y_0<\ell_{\pi_{k+1}}$} \Comment{Feasibility Phase} \label{line:start of feasibility phase}
		\State $(y_0,\bm\alpha,\bm\beta,\gamma)\gets\pivot{y_0}{\ell_{\pi_{k+1}}}$
		\EndWhile \label{line:end of feasibility phase}
		\If{$\nabla_{\pi_{k+1}}f\left(\bm{\ell}_{\underline{\alpha}},\nabla_\gamma^{-1} f(\bm a_{\gamma};\bm{\ell}_{\underline{\alpha}},\bm u_{\overline{\alpha}},\bm 0_{\gamma_0},y_0),\bm u_{\overline{\alpha}},\bm 0_{\gamma_0},\ell_{\pi_{k+1}}\right)- a_{i}\ge0$} \label{line:start of optimality phase}
		\State $\alpha_-\gets\alpha_-\cup\{ \pi_{k+1} \}$  \Comment{Optimality Phase}
		\Else
		\Repeat
		\State $\bar y\gets$  root of $\nabla_{\pi_{k+1}}f\left(\bm{\ell}_{\underline{\alpha}},\nabla_\gamma^{-1} f(\bm a_{\gamma};\bm \ell_{\underline{\alpha}},\bm u_{\overline{\alpha}},\bm 0_{\gamma_0},y),\bm u_{\overline{\alpha}},\bm 0_{\gamma_0},y\right)- a_{i}=0$ \label{line:root1}
		\State $ \bar y\gets \min\{ \bar y, u_{\pi_{k+1}} \}$ \label{line:maintain feasibility}
		\State  $(y_0,\bm\alpha,\bm\beta,\gamma)\gets\pivot{y_0}{\bar y}$ \label{line:trace path}
		\Until{$y_0=\bar y$} \Comment{$y=\bar y$ is optimal for \eqref{eq:subproblem-y}}
		\If{$\bar y<u_{\pi_{k+1}}$} 
		\State $\alpha_+\gets\alpha_+\cup\{\pi_{k+1}\}$
		\EndIf
		\If{$\bar y=u_{\pi_{k+1}}$}
		\State $\overline{\alpha}\gets\overline{\alpha}\cup\{ \pi_{k+1} \}$
		\EndIf 
		\EndIf \label{line:end of optimality pahse}
		\EndIf
		\If{$n_1<\pi_{k+1}\le n$}
		
		\If{$y_0=0$}
		\State $\beta_\oplus\gets\beta_\oplus\backslash\{\pi_{k+1}\}$
		\Else \Comment{$y_0<0$} \label{line:start of feasibility case 2}
		\State $\beta_-\gets\beta_-\backslash\{ \pi_{k+1} \}$
		\Repeat 
		\State $(y_0,\bm\alpha,\bm\beta,\gamma)\gets\pivot{y_0}{0}$ 
		\Until{$y_0=0$} 
		\EndIf		\label{line:end of feasibility case 2}
		\State $\beta_0\gets\beta_0\cup\{\pi_{k+1}\}$
		\EndIf		

		\If{$n<\pi_{k+1}\le 2n-n_1$} \Comment{$y_0=0$}
		\State $\beta_0\gets\beta_0\backslash\{ \pi_{k+1} \}$
			\If{$\nabla_{\pi_{k+1}}f\left(\bm{\ell}_{\underline{\alpha}},\nabla_\gamma^{-1} f(\bm a_{\gamma};\bm \ell_{\underline{\alpha}},\bm u_{\overline{\alpha}},\bm 0_{\gamma_0},y),\bm u_{\overline{\alpha}},\bm 0_{\gamma_0},0_{\pi_{k+1}}\right)- a_{i}\ge0$}
		\State $\beta_\oplus\gets\beta_\oplus\cup\{ \pi_{k+1} \}$ \Comment{$y=0$ is optimal for \eqref{eq:subproblem-y} }
		\Else
		\Repeat 
		\State $\bar y\gets$  root of $\nabla_{\pi_{k+1}}f\left(\bm{\ell}_{\underline{\alpha}},\nabla_\gamma^{-1} f(\bm a_{\gamma};\bm \ell_{\underline{\alpha}},\bm u_{\overline{\alpha}},\bm 0_{\gamma_0},y),\bm u_{\overline{\alpha}},\bm 0_{\gamma_0},y\right)- a_{i}=0$ \label{line:root2}
		\State  $(y_0,\bm\alpha,\bm\beta,\gamma)\gets\pivot{y_0}{\bar y}$
		\Until{$y_0=\bar y$} \Comment{$y=\bar y$ is optimal for \eqref{eq:subproblem-y}}
		\State $\beta_+\gets\beta_+\cup\{\pi_{k+1}\}$
		\EndIf
		\EndIf
		\State $\gamma\gets\alpha_+\cup\beta_-\cup\beta_+$, $x\gets$solution defined in \eqref{eq:index->x}, $v(\ones^{\pi_{[k+1]}})\gets f(x)-c^\top x$ \label{line:solution}
		\EndFor
		\State \Return $\left\{ v(\ones^{\pi_{[k]}}) \right\}_{k=0}^{2n-n_1}$
	\end{algorithmic}
\end{algorithm}

Proposition~\eqref{prop:linear-step-algo} shows that Algorithm~\ref{algo:active-set} can correctly compute the extreme bases and encounters $\bigO{n}$ breakpoints during its execution.
\begin{proposition}\label{prop:linear-step-algo}
	Algorithm~\ref{algo:active-set} computes the extreme base  $\left\{ v(\ones^{\pi_{[k]}}) \right\}_{k=0}^{2n-n_1}$ with $\bigO{n}$ invocations of $\pivot{y_0}{\bar y}$.
\end{proposition}
\begin{proof}
	Define $\tilde{x}$ as the solution obtained from Line~\ref{line:solution} in Algorithm~\ref{algo:active-set}. Because throughout the implementation of each for-loop $\bm\alpha$ and $\bm\beta$ always give rise to the optimal solution to \eqref{eq:subproblem-y}, according to Lemma~\ref{lem:isotonicity} and in order to prove the correctness of the algorithm at termination, it suffices to show that $\tilde x_{\pi_{k+1}}=\argmin_{y\in\mathcal{I}_{k+1}}v_{k+1}^\pi(y)$, where  $\mathcal{I}_{k+1}=[\ell_{\pi_{k+1}},u_{\pi_{k+1}}]$ if $\pi_{k+1}\in[n_1]$, $\mathcal{I}_{k+1}=\{0\}$ if $\pi_{k+1}\in[n]\backslash[n_1]$, and $\mathcal{I}_{k+1}=[0,\infty)$ if $\pi_{k+1}\in[2n-n_1]\backslash[n]$, corresponding to the three major if-cases in the $k$-th for-loop. Next, we perform a case-by-case analysis. 
	\begin{itemize}[leftmargin=*]
	\item \emph{Case 1:}  $\pi_{k+1}\in[n_1]$. In this case, $\pi_{k+1}\in\mathcal{N_+}$, $z_{\pi_{k+1}}$ is one of $z^+$-variables,  $\bar x^k_{\pi_{k+1}}=0$, and $\mathcal{I}_{k+1}=[\ell_{\pi_{k+1}},u_{\pi_{k+1}}]$. Because $\ell_{\pi_{k+1}}\ge0$, in the beginning of Line~\ref{line:start of the loop} one has $y_0=0\not\ge\ell_{\pi_{k+1}}$ unless $\ell_{\pi_{k+1}}=0$. Therefore, the algorithm consists of two phases -- Feasibility Phase (Line~\ref{line:start of feasibility phase}--Line~\ref{line:end of feasibility phase}) and Optimality Phase (Line~\ref{line:start of optimality phase}--Line~\ref{line:end of optimality pahse}). In the Feasibility Phase, one increases $y_0$ and calls $\pivot{y_0}{\bar y}$ to trace the solution path of until $y_0=\ell_{k+1}$.  At this point, $y_0$ becomes feasible and thus, we turn to the Optimality Phase to seek the optimal $y$ over $\mathcal{I}_{k+1}$. 
	
	Because $v^\pi_{k+1}(\cdot)$ is a convex function and $\left(v^\pi_{k+1}\right)'(y)=\nabla_{\pi_{k+1}} f(\bm x^{k+1}(y))$, the optimality condition of minimizing the value function $\min\{ v_{k+1}^\pi(y): \ell_{\pi_{k+1}}\le y \le u_{\pi_{k+1}}\}$ is given by 
	\[ \begin{cases}
		y=\ell_{\pi_{k+1}}&\text{ if } \nabla_{\pi_{k+1}} f(x^{k+1}(y))\ge0 \\
		y =\bar y\text{ satisfying }\nabla_{\pi_{k+1}} f(\bm x^{k+1}(y))=0 &\text{ if } \ell_{\pi_{k+1}}<y< u_{\pi_{k+1}}\\
		y=u_{\pi_{k+1}} &\text{ otherwise.}
	\end{cases} \]
	Note that $(\bm x(y))_\gamma=\nabla_\gamma^{-1} f(\bm a_{\gamma};\bm \ell_{\underline{\alpha}},\bm u_{\overline{\alpha}},\bm 0_{\gamma_0},y)$. Thus, if the condition in Line~\ref{line:start of optimality phase} fails, then one can deduce that the optimal $y>\ell_{\pi_{k+1}}$ and we increase $y$. When $y$ runs over $(\ell_{\pi_{k+1}},u_{\pi_{k+1}})$, we keep tracing the solution path (Line~\ref{line:trace path}) and meanwhile maintain $y\le u_{\pi_{k+1}}$ (Line~\ref{line:maintain feasibility}) until $y=\bar y$ or otherwise, we must have the optimal $y=u_{\pi_{k+1}}$. This proves the correctness of the algorithm in the first case.
	\item \emph{Case 2: $\pi_{k+1}\in[n]\backslash[n_1]$.} In this case, $\pi_{k+1}\in\mathcal{N_\pm}$ and $z_{\pi_{k+1}}$ is one of $z^-$-variables. Because $\mathcal{I}_{k+1}$ is a singleton,  only the Feasibility Phase (Line~\ref{line:start of feasibility case 2}--Line~\ref{line:end of feasibility case 2}) is required. 
	\item \emph{Case 3: $\pi_{k+1}\in[2n-n_1]\backslash[n]$}. In this case,  $\pi_{k+1}-n+n_1\in\mathcal{N_\pm}$ and $z_{\pi_{k+1}}$ is one of $z^+$-variables. Also, one has $\mathcal{I}_{k+1}=[0,\infty)$ and $y_0=0$, implying that $y_0$ is already feasible and only the Optimality Phase is required. Because the analysis in Case~2 and Case~3 is similar to that of Case~1, the details are omitted.
	\end{itemize}
	
	Finally, we prove the linear complexity of the algorithm. This follows from that throughout Algorithm~\ref{algo:active-set}, the state of the variable $x_i$ for $i\in[n_1]$ can only transit along the path  $\alpha_0\to\underline{\alpha}\to\alpha_+\to\overline{\alpha}$. Additionally, the state of the variable $x_i$ for $i\in[n]\backslash[n_1]$ can only transit along the path $\beta_-\to\beta_\ominus\to\beta_0\to\beta_\oplus\to\beta_+$. Some edges in the two paths might be skipped. Consequently, the transition can occur at most $3n_1+4(n-n_1)=4n-n_1=\bigO{n}$ times. This finishes the proof.
\end{proof}

In general, the implementation of Algorithm~\ref{algo:active-set} relies on computing $\nabla_\gamma ^{-1}f(\cdot)$ and solving a univariate root-finding problem. The former amounts to solving an unconstrained convex program, for which plenty of convex optimization algorithms are applicable. Furthermore, since $f$ is submodular, specialized algorithms (\cite{chandrasekaran1970special}, \cite{tamir1974minimality})also exist that can further improve computational efficiency. Regarding the latter, since all univariate equations arising in Algorithm~\ref{algo:pivot} and Algorithm~\ref{algo:active-set} are monotonic, their roots can be found easily by employing standard numerical methods. In some special cases, such as when $f$ is quadratic, these related quantities admit an analytical form. In the sequel, we focus on the specialization of the algorithm to quadratic and conic quadratic $f$.
\subsection{Tracing solutions paths in quadratic cases.}
Assume $f(\bm x)=\frac{1}{2}\bm x^\top \bm Q\bm x$, where $\bm Q$ is a Stieltjes matrix. In this scenario, Algorithm~\ref{algo:pivot} is closely tied to pivoting methods for solving linear complementarity problems. For example,  Line~\ref{line:ratio test} is an analogy of the ratio-test operation; see \cite{cottle2009linear}, Chapter~4 for more  details. 

   We now present the closed-form expressions of key quantities in Algorithm~\ref{algo:pivot} and Algorithm~\ref{algo:active-set}. First, one can verify that $\bm x_\gamma$ in \eqref{eq:index->x} is given by
   \[ \bm x_\gamma=\nabla_\gamma^{-1}f\left( \bm a_\gamma;\bm{\ell}_{\underline{\alpha}},\bm u_{\overline{\alpha}},\bm 0_{\gamma_0},y \right)=-\bm Q_{\gamma\gamma}^{-1}\bm Q_{\gamma\pi_{k+1}}y+\bm Q^{-1}_{\gamma\gamma}(\bm a_\gamma-\bm Q_{\gamma\underline{\alpha}}\bm{\ell}_{\underline{\alpha}}-\bm Q_{\gamma\overline{\alpha}}\bm u_{\overline{\alpha}}).  \]
   With the expression of $\nabla_\gamma^{-1}f\left( \bm a_\gamma;\bm{\ell}_{\underline{\alpha}},\bm u_{\overline{\alpha}},\bm 0_{\gamma_0},y \right)$, the ratios $r_i$  in Algorithm~\ref{algo:pivot} can be readily calculated
\begin{align*}
	r_i=\begin{cases}
		\frac{\bm{Q}^{-1}_{\gamma\gamma}(\bm{a}_\gamma-\bm{Q}_{\gamma\underline{\alpha}}\bm{\ell}_{\underline{\alpha}}-\bm{Q}_{\gamma\overline{\alpha}}\bm u_{\overline{\alpha}})-u_i}{\bm{Q}_{\gamma\gamma}^{-1}\bm{Q}_{\gamma\pi_{k+1}}}&\text{if }i\in\alpha_+\\
		\frac{\bm{Q}^{-1}_{\gamma\gamma}(\bm{a}_\gamma-\bm{Q}_{\gamma\underline{\alpha}}\bm{\ell}_{\underline{\alpha}}-\bm{Q}_{\gamma\overline{\alpha}}\bm u_{\overline{\alpha}})}{\bm{Q}_{\gamma\gamma}^{-1}\bm{Q}_{\gamma\pi_{k+1}}}&\text{if }i\in\beta_-\\
		\frac{a_i-\bm{Q}_{i\gamma}\bm{Q}_{\gamma\gamma}^{-1}\bm{a}_\gamma-\left( \bm{Q}_{i\underline{\alpha}}-\bm{Q}_{i\gamma}\bm{Q}_{\gamma\gamma}^{-1}\bm{Q}_{\gamma\underline{\alpha}} \right)\bm{\ell}_{\underline{\alpha}}-\left( \bm{Q}_{i\overline{\alpha}}-\bm{Q}_{i\gamma}\bm{Q}_{\gamma\gamma}^{-1}\bm{Q}_{\gamma\overline{\alpha}} \right)\bm u_{\overline{\alpha}}}{\bm{Q}_{i\pi_{k+1}}-\bm{Q}_{i\gamma}\bm{Q}_{\gamma\gamma}^{-1}\bm{Q}_{\gamma\pi_{k+1}}}&\text{if }i\in\underline{\alpha}\cup\beta_\oplus.
	\end{cases}
\end{align*}
 Additionally, the root of equations in Line~\ref{line:root1} and Line~\ref{line:root2} of Algorithm~\ref{algo:active-set} share the same formula
\[  \bar{y}= 	\tfrac{a_{\pi_{k+1}}-\bm{Q}_{{\pi_{k+1}}\gamma}\bm{Q}_{\gamma\gamma}^{-1}\bm{a}_\gamma-\left( \bm{Q}_{{\pi_{k+1}}\underline{\alpha}}-\bm{Q}_{{\pi_{k+1}}\gamma}\bm{Q}_{\gamma\gamma}^{-1}\bm{Q}_{\gamma\underline{\alpha}} \right)\bm{\ell}_{\underline{\alpha}}-\left( \bm{Q}_{{\pi_{k+1}}\overline{\alpha}}-\bm{Q}_{{\pi_{k+1}}\gamma}\bm{Q}_{\gamma\gamma}^{-1}\bm{Q}_{\gamma\overline{\alpha}} \right)\bm u_{\overline{\alpha}}}{\bm{Q}_{{\pi_{k+1}}\pi_{k+1}}-\bm{Q}_{{\pi_{k+1}}\gamma}\bm{Q}_{\gamma\gamma}^{-1}\bm{Q}_{\gamma\pi_{k+1}}}.\]
Here, each  $r_i$ or $\bar y$ is understood as $\infty$ if the denominator of the ratio is 0. 

 Note that in each iteration the state of only one index is changed, leading to a rank-one update of $\bm{Q}_{\gamma\gamma}^{-1}$. Consequently, the computation of key quantities listed above can be accomplished in $\bigO{n^2}$ time per iteration in an incremental way. We refer readers to  \cite{gomez2022linear, he2023comparing, pang2023some} and the references therein for details.  In addition, one can solve the initial subproblem to get $v(\bm 0)$ and the associated optimal solution in $\bigO{n^3}$ \cite{pang2023some}. Combining this fact with the quadratic complexity per step and the linear number of steps (Proposition~\ref{prop:linear-step-algo}), one obtains the overall complexity of computing the extreme bases in the quadratic case.
\begin{proposition}
	If $f(x)=\frac{1}{2}\bm x^\top \bm{Q}\bm x$ with $\bm{Q}$ being a Stieltjes matrix, Algorithm~\ref{algo:active-set} can terminate in $\bigO{n^3}$ time.
\end{proposition}

Notably, in this case, the cubic complexity matches the best known complexity of computing $v(\one)$, thus the extreme bases can be computed in the same complexity as an evaluation of the continuous submodular function. We conclude this section with an example to illustrate Algorithm~\ref{algo:active-set}.
\begin{example}
	Consider 
	\[  \bm{Q}=\begin{bmatrix}
		5&-1&-3\\
		-1&3&-2\\
		-3&-2&7
	\end{bmatrix},\;\bm a=\begin{bmatrix}
	1\\1\\1
	\end{bmatrix},\; \bm\ell=\begin{bmatrix}
	0\\0\\0
	\end{bmatrix},\; \bm u=\begin{bmatrix}
	1\\1\\1
	\end{bmatrix}. \]
	Given above data and permutation $\pi=(1,2,3,4)$,  one can compute the extreme bases $\left\{v(0,0,0),v(1,0,0),v(1,1,0),v(1,1,1)\right\}$ using Algorithm~\ref{algo:active-set}. It can be seen easily that $v(0,0,0,0)=0$. The solution $(x_1(y),x_2(y),x_3(y))$ to subproblems $P^1(y),P^2(y) $ and $P^3(y)$ are shown in Figure~\ref{fig:trajectory}. For each $k=1,2,3$, we use $y\gets x_k$ as the driving parameter to drive the increase of the solution $\bm x(y)$. In this example, one encounters four breakpoints during the implementation of Algorithm~\ref{algo:active-set}.
	\begin{figure}[H]
		\includegraphics[scale=0.3]{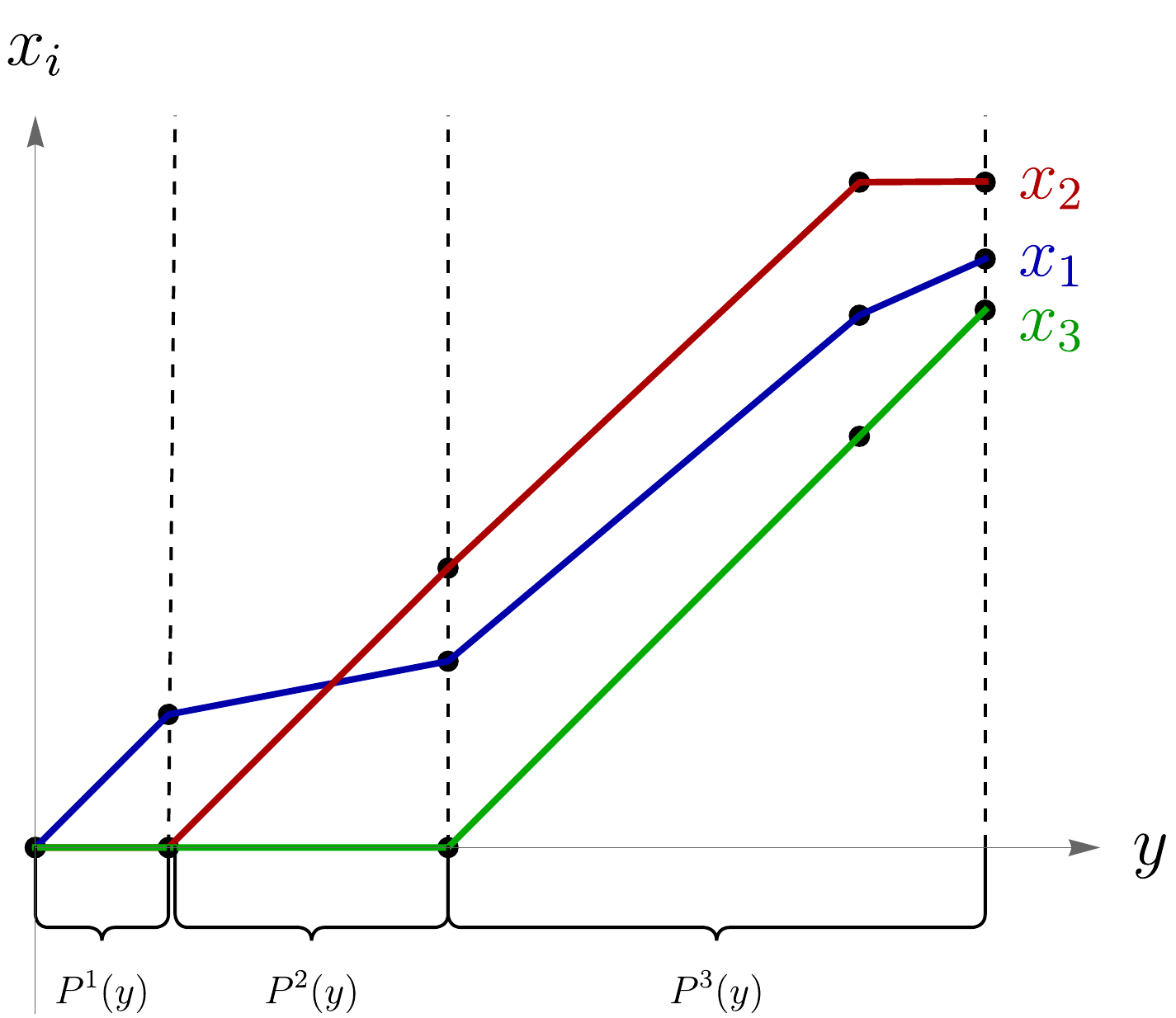}
		\caption{Trajectory of $x_i$ induced by Algorithm~\ref{algo:active-set}}\label{fig:trajectory}
	\end{figure}
\end{example}

\subsection{Tracing solutions paths in conic diagonal quadratic cases.} In this section, we consider the problem of minimizing
\begin{equation}\label{eq:cqp-diagonal}
	\begin{aligned}
		\min_{x,z}\;&\sqrt{\sigma^2+\sum_{i=1}^nc_ix_i^2}-\bm a^\top \bm x+ \bm d^\top \bm z\\
		\text{s.t. }&0\le x_i\le u_iz_i\;\forall i=1,\dots,n,
	\end{aligned} 
\end{equation}
where $f(x)=\sqrt{\sigma^2+\sum_{i=1}^nc_ix_i^2}$ is a convex submodular function over $\R_+^n$, $\sigma>0$, and $\bm c,\bm u>\bm 0$. Since $a_i\le0$ implies $x_i=0$ at optimality, we assume $\bm a>\bm 0$ without loss of generality. By rescaling each $x_i$, we may further assume $c_i=1$ for all $i\in[n]$. Under this setup, we have $\gamma=\alpha_+$ and all index sets $\beta$'s vanish. Moreover, $\alpha_0$ is unnecessarily needed as $\bm\ell=\bm 0$. 

We now specify the key quantities in the solution path tracing algorithm. Because $\nabla_\gamma f(\bm x)=\bm x_\gamma/f(x)$, we obtain from \eqref{eq:index->x}
\[ \bm x_\gamma=\nabla_\gamma^{-1}f\left( \bm a_\gamma;\bm 0_{\underline{\alpha}},\bm u_{\overline{\alpha}},y \right)=\sqrt{\frac{\sigma^2+\norm{2}{\bm u_{\overline{\alpha}}}^2+y^2}{1-\norm{2}{\bm a_\gamma}^2}}\bm a_\gamma. \]
Because $a_i>0$ and $\ell_i=0$ implies $\nabla_i f(\bm\ell_{\underline{\alpha}},\bm x_\gamma,\bm u_{\overline{\alpha}},y)$=0 for all $i\in\underline{\alpha}$, the equation in Line~\ref{line:lb-ri} of Algorithm~1 yields no root. Hence, we only calculate $r_i$ for $i\in\alpha_+$ which is given by $r_i=h(u_i/a_i)$, where 
\[ h(p)\defeq\sqrt{\left(1-\norm{2}{\bm a_\gamma}^2\right)p^2-\norm{2}{\bm u_{\overline{\alpha}}}^2-\sigma^2}. \]
 Similarly, the root in Line~\ref{line:root1} of Algorithm~\ref{algo:active-set} is given by $\bar y=h(u_{\pi_{k+1}}/a_{\pi_{k+1}})$. 

Above simplifications lead to a streamed implementation of Algorithm~\ref{algo:pivot} and \ref{algo:active-set}. The next breakpoint in this case is computed as
\begin{equation}\label{eq:bp-cqp}
	 y_{\texttt{bp}}\gets\min\left\{ u_{\pi_{k+1}}, h\left(\frac{u_{\pi_{k+1}}}{a_{\pi_{k+1}}}\right),\min_{i\in\alpha_+}h\left(\frac{u_i}{a_i}\right) \right\}. 
\end{equation}

Although we assume $\sigma>0$ to ensure differentiability of $f(\cdot)$, the algorithm remains valid for $\sigma=0$. Note that each pivoting operation costs $\bigO{n\log(n)}$ due to sorting over at most $n$ values. Combined with Proposition~\ref{prop:linear-step-algo}, the total complexity of Algorithm~\ref{algo:active-set} is $\bigO{n^2\log(n)}$ under a naive implementation in this scenario. Proposition~\ref{prop:cqp} shows that with suitable data structures, the time complexity can be reduced by a factor of $n$.
\begin{proposition}\label{prop:cqp}
	If $f(x)=\sqrt{\sigma^2+\sum_{i=1}^nc_ix_i^2}$ with $\bm c>0$ and $\sigma\ge0$, then Algorithm~\ref{algo:active-set} terminates in $\bigO{n\log(n)}$ time.
\end{proposition}
\begin{proof}
	Because $h(\cdot)$ is an increasing function over $[0,\infty)$, computing the next breakpoint from \eqref{eq:bp-cqp} reduces to comparing $u_i/a_i$ for $i\in\alpha_+\cup\{\pi_{k+1}\}$. This can be accomplished by maintaining a \emph{priority queue} that stores the sorted sequence of $\{ u_i/a_i \}_{i\in\alpha_+}$. When encountering a new breakpoint, the potential update to $\alpha_+$ modifies the priority queue via a single insertion or deletion, both taking $\bigO{\log(n)}$ time. The result follows from Proposition~\ref{prop:linear-step-algo}.
\end{proof}
\citet{atamturk2019lifted} propose an algorithm for solving $v(\bm 1)$  in $\bigO{n\log(n)}$ in the conic diagonal quadratic setting. In contrast, Algorithm~\ref{algo:active-set} is able to solve all $n$ subproblems with the same time complexity.

%

\section{Computations}\label{sec:computations}
In this section, we evaluate the performance of different exact global solution methods for tackling combinatorial quadratic MRF inference problems. Our focus narrows down to two distinct classes of exact-solution approaches: binary submodular minimization and mixed-integer programming (MIP). The following subsections delve into our investigation. Section~\ref{sec:mrf-nonnegative} is devoted to synthetic instances of sparse MRF problems, assuming nonnegative lower bounds on the continuous variables.  In Section~\ref{sec:mrf-free}, we shift our attention to an outlier detection problem using time series data from the CBLIB library. This segment involves a comparative study between the methodology developed in this paper and the state-of-the-art MIP approach.

\subsection{Sparse MRF inference}\label{sec:mrf-nonnegative}
Consider a general sparse quadratic program with a Stieltjes matrix 
\begin{equation}\label{eq:mrf-nonnegative}
	\begin{aligned}
		\min_{\bm x,\bm z}\;&\frac{1}{2}\bm{x}^\top \bm{Q} \bm{x}-\bm{a}^\top \bm{x}+\bm{c}^\top \bm{z}\\
		\text{s.t. }\;&\bm{\ell}\circ \bm{z}\le \bm{x}\le\bm{u}\circ\bm{z},\;\bm{z}\in\{0,1\}^n,
	\end{aligned}
\end{equation}
where $\bm{0}<\bm{\ell}<\bm{u}$ are $n$-dimensional vectors, $\bm Q\in\R^{n\times n}$ is a Stieltjes matrix, and $\bm{a},\bm{c}\in\R^n$. It can be shown that for any Stieltjes matrix $Q$, there always exist quadratic functions $h_{ij}$ and $g_i$ such that $x^\top Q x$ can be decomposed as the sum of one/two-dimensional forms. For this reason, \eqref{eq:mrf-nonnegative} can be put in the form of \eqref{eq:SparseMRF}, and we can interpret \eqref{eq:mrf-nonnegative} as a sparse MRF inference problem. 
 All experiments in this section are conducted on a laptop with a 2.30GHz $\text{Intel}^\text{\textregistered}$
 $\text{Core}^{\text{\tiny TM}}$ i9-9880H CPU and 64 GB main memory. 
 \subsubsection{Instance generation}\label{sec:instance-generation}
 We now describe how we generate synthetic instances. Given dimension $n$, the data tuple $(\bm Q,\bm a,\bm c,\bm \ell,\bm u)$ is generated in the following way:
 \begin{itemize}
 	\item Draw $n$ numbers $\hat c_i$ independently from normal distribution $\set N(0,7e5)$. Let $c_i=|\hat c_i|$ for all $i\in[n]$.
 	\item Draw $n$ numbers $\hat a_i$ independently from normal distribution $\set N(0,1e5)$. Let $a_i=|\hat a_i|$ for all $i\in[n]$.
 	\item Set $\ell_i=2$ and $u_i=10$ for all $i\in[n]$.
 	\item For each $i\in[n]$ and $j\in[n]$, draw $W_{ij}$ independently and uniformly from $[0,1]$. Let $M_{ij}=-|W_{ij}+W_{ji}|/2$ for all $i\in[n]$, $j\in[n]$. Set $Q_{ij}=M_{ij}$ for $i\neq j$ and $Q_{ii}=\sum_{j\neq i}|M_{ij}|$ for $i\in[n]$ to ensure that $\bm{Q}$ is a Stieltjes matrix. 
 \end{itemize}
\subsubsection{Efficiency of computing extreme bases}\label{sec:computation-lovasz}
Computation of Lov\'{a}sz extension plays a pivotal role in solving \eqref{eq:mrf-nonnegative} as a binary submodular minimization problem. This section is dedicated to the evaluation of different methods employed to compute the extreme bases associated with \eqref{eq:mrf-nonnegative}. Specifically, for each $n\in\{50,80,100,200,300,400,500,1000\}$, we generate five instances as outlined in Section~\ref{sec:instance-generation}. For each instance generated, we compute the extreme bases $\{v_k\}_{k\in[n]}$, where \[v_k\defeq\min\left\{\dfrac{1}{2}\bm{x}^\top \bm{Q} \bm{x}-\bm{a}^\top \bm{x}:\begin{aligned}
	&\ell_i\le x_i\le u_i\,\forall i\in[k],\\
	&x_i=0\,\forall i\in[n]\backslash[k] 
\end{aligned}\right\},\] using the following three methods 
\begin{itemize}
	\item \slow: We directly solve a series of $n$ convex quadratic programs that define $v_k$ using Gurobi;
	\item \fastp: We use Algorithm~\ref{algo:active-set} to progressively compute $\{v_k\}_{k\in[n]}$.
\end{itemize}

\begin{table}[h]
	\caption{Time for computation of extreme bases}\label{tab:computing-lovasz}
	\centering
	\begin{tabular}{c|ccc}
		\toprule
		dimension&\slow&\fastp\\\midrule
		50& 0.095& 0.0008\\
		80&0.2&0.0016\\
		100&0.36&0.0045\\
		200&2.59&0.04\\
		300&7.94&0.12\\
		400&22.75&0.56\\
		500&44.92&1.2\\
		1000&379.31&13.67\\\bottomrule
	\end{tabular}
\end{table}
Each row of Table~\ref{tab:computing-lovasz} presents the average computation time, measured in seconds, for obtaining $\{v_k\}_{k\in[n]}$ over five instances. Notably,  the performance of \slow stands out as the least efficient, exhibiting a 30-fold increase in execution time compared to \fastp across all scenarios.  
\subsubsection{Results in solving MRFs}\label{sec:mrf-nonnegative-results}
We test the performance of the following two methods to solve \eqref{eq:mrf-nonnegative}
\begin{itemize}
	\item \mip: solve \eqref{eq:mrf-nonnegative} as a mixed-integer program using Gurobi~9.0.2 with default settings;
	\item \sm(
	\slow): solve \eqref{eq:mrf-nonnegative} as a binary submodular minimization problem using the cutting plane method which is included in Appendix~\ref{sec:cutting-plane} for completeness.  Lov\'asz extensions are computed by solving quadratic programs using Gurobi.
	\item \sm(\fastp): solve \eqref{eq:mrf-nonnegative} as a binary submodular minimization problem using the cutting plane method, using Algorithm~\ref{algo:active-set} as a subroutine.
\end{itemize}

The time limit for both solution methods is set as 1800 seconds. Each entry of Table~\ref{tab:positive_MRF} shows average statistics over five instances. The notations \slow and \fastp signify the methodologies employed for computing the extreme bases associated with  \eqref{eq:mrf-nonnegative}. The table displays the dimension of the problem $n$, the solution time for solving \eqref{eq:mrf-nonnegative} (Time), the final gap reported by the solver upon termination (Gap), the count of instances solved to optimality within the prescribed time limit (\#), the proportion of solving time attributed to the computation extreme bases (EB), and the sparsity of optimal solutions quantified as $\text{\sparsity}=\dfrac{\sum_{i=1}^n z^*_i}{n}\times 100\%$, where $z_i^*$ represents the optimal indicator variables pertaining to \eqref{eq:mrf-nonnegative}.

In Table~\ref{tab:positive_MRF}, it is evident that \mip exhibits the poorest performance when considering the number of instances solved to optimality. It can solve all instances with dimension $n\le 200$ but none of high-dimensional instances. \sm(\slow) is marginally better, resolving one additional instance with $n=300$. However, \sm(\slow) still struggles with most high-dimensional ones, and its average solution time of \sm(\slow) lags by at least a factor of seven compared to \mip.  In stark contrast, \sm(\fastp) excels in both solvability and solution time, managing to tackle all test instances within half of the allotted time limit.  For the instances solvable by both \mip and \sm(\slow), \sm(\fastp) is able to solve them in mere five seconds, signifying a remarkable improvement. Because the only distinction between \sm(\fastp) and \sm(\slow) lies in the way of evaluating extreme bases, one can conclude that Algorithm~\ref{algo:active-set} plays a significant role in the success of \sm(\fastp). This is also consistent with the conclusion from Section~\ref{sec:computation-lovasz} and the observation that \sm(\fastp) typically spends over 95\% running time on the computation of the Lov{\'a}sz extension. 
\begin{table}[h]
	\centering
	\caption{Results for solving sparse MRF inference problems}\label{tab:positive_MRF}
	\resizebox{\linewidth}{!}{
		\begin{tabular}{r|rcc|rcc|rcc|l}
			\toprule		
			\multirow{2}{*}{$n$} &\multicolumn{3}{c|}{\mip}&\multicolumn{3}{c|}{\sm(\slow)}&\multicolumn{3}{c|}{\sm(\fastp)}&\multirow{2}{*}{\sparsity} \\\cmidrule(lr){2-4} \cmidrule(lr){5-7} \cmidrule(lr){8-10}
			& Time(s) & Gap(\%) & \# &   Time(s) &  EB(\%) & \# &  Time(s) &  EB(\%) & \#  & \\\midrule
			
			50 &       0.02 &       0 &          5 &       1.54 &      99.53 &          5 &       0.01 &      87.91 &          5 &       27.2 \\
			
			80 &       0.04 &       0 &          5 &       3.99 &      99.77 &          5 &       0.02 &      90.27 &          5 &      27.75 \\
			
			100 &       0.06 &       0 &          5 &       8.91 &      99.92 &          5 &       0.08 &      96.42 &          5 &       20.6 \\
			
			200 &      40.83 &       0 &          5 &     285.95 &      99.96 &          5 &       4.32 &      97.53 &          5 &       18.1 \\\midrule
			
			300 &    1800.00 &    11.25 &          0 &    1800.00 &      99.97 &          1 &      52.00 &      97.14 &          5 &       14.2 \\
			
			400 &    1800.00 &    34.89 &          0 &    1800.00 &      99.99 &          0 &     197.74 &      98.37 &          5 &      10.05 \\
			
			500 &    1800.00 &    52.81 &          0 &    1800.00 &     100.00 &          0 &     868.57 &      98.47 &          5 &       8.52 \\\bottomrule
			
		\end{tabular}
	}	
\end{table}
\subsection{Outlier detection in time series}\label{sec:mrf-free}
Given data $(\bm{\tau,\mu,y,\sigma})\in\R^{n\times 4}$ with time stamps $0<\tau_1<\dots<\tau_n$, consider the problem of outlier detection in time series of the form 
\begin{equation}\label{eq:outlier-detection}
	\hspace{-1em}\begin{aligned}		
		\min_{\bm{x},\bm{w},\bm{z}}\;&\frac{x_1^2}{2\tau_1}+\sum_{i=1}^{n-1}\frac{\left(x_{i+1}-x_i\right)^2}{2(\tau_{i+1}-\tau_i)}+\sum_{i=1}^n\frac{(y_i+w_i-\mu_i-x_i)^2}{2\sigma_i^2}+c\sum_{i=1}^nz_i\\
		\text{s.t. }\;&\bm{w}\circ(\ones-\bm{z})=0,\;\bm{z}\in\{0,1\}^n,
	\end{aligned}
\end{equation}
where $c$ is the parameter controlling the number of outliers to be discarded. Note that \eqref{eq:outlier-detection} is a special case of robust MRF inference problems introduced in   Section~\ref{sec:robust-mrf-inference}. For the background and statistical model regarding \eqref{eq:outlier-detection}, we refer readers to \cite{gomez2021outlier}.

\subsubsection{Solution methods}
We now outline the three solution methods employed in this section to tackle \eqref{eq:outlier-detection}. The first method corresponds to using standard big-M formulation of the problem with a MIO solver, and the second consists of using a strong conic formulation proposed in \cite{gomez2021outlier}, which represents the current state-of-the-art MIO formulation. 
The third method is the submodular minimization approach introduced in this work. 

\noindent$\bullet$ \bigM. In \bigM, we reformulate the complementarity constrains using standard big-M techniques with $\displaystyle M\defeq \max_{j\in[n]}\{y_j-\mu_j\}-\min_{j\in[n]}\{y_j-\mu_j\}$
\begin{equation}
	\hspace{-1em}\begin{aligned}		
		\min_{\bm{x},\bm{w},\bm{z}}\;&\frac{x_1^2}{2\tau_1}+\sum_{i=1}^{n-1}\frac{\left(x_{i+1}-x_i\right)^2}{2(\tau_{i+1}-\tau_i)}+\sum_{i=1}^n\frac{(y_i+w_i-\mu_i-x_i)^2}{2\sigma_i^2}+c\sum_{i=1}^nz_i\\
		\text{s.t. }\;&-M\bm{z}\le\bm{w}\le M\bm{z},\;\bm{z}\in\{0,1\}^n.
	\end{aligned}
\end{equation}

\noindent$\bullet$ \strongMip.  In \strongMip, we adopt the strong mixed-integer formulation of \eqref{eq:outlier-detection} based on convexification techniques
\begin{equation}\label{eq:strong mip}
	\hspace{-1em}\begin{aligned}
		\min_{\bm{x,z,w,s,\bar z, r}}\;&\frac{x_1^2}{2\tau_1}+\frac{1}{2}\sum_{i=1}^{n-1}\left( \lambda_is_{i,1}^2+\frac{(s_{i,1}-s_{i,2})^2}{\tau_{i+1}-\tau_i}\right.\\
		&\left.+\left( \frac{1}{\sigma_{i+1}^2}-\lambda_{i+1} \right)s_{i,2}^2+\frac{\lambda_i(1/\sigma_{i+1}^2-\lambda_{i+1})}{L_i}r_i \right)\\
		&-\sum_{i=1}^n\frac{(y_i-\mu_i)(x_i-w_i)}{\sigma_i^2}+c\sum_{i=1}^nz_i + \sum_{i=1}^n\frac{(y_i-\mu_i)^2}{2\sigma_i^2}\hspace{-2em}\\
		\text{s.t. }\;& s_{i,1}=x_i-w_i+\frac{1/\sigma_{i+1}^2-\lambda_{i+1}}{L_i}(2w_i-w_{i+1}),& i\in[n-1] \\
		&s_{i,2}=x_{i+1}-v_{i+1}-\frac{\lambda_i}{L_i}(w_i-w_{i+1}),&i\in[n-1]\\
		&\bar z_i\le 1,\;\bar{z}_i\le z_i+z_{i+1},\;(w_i-w_{i+1})^2\le r_i\bar{z}_i,& i\in[n-1]\\
		&-M\bm{z}\le \bm{w}\le M\bm{z},\;\bm{z}\in\{0,1\}^n,\;\bm{w}\in\R^n,\\
		&\bm{x}\in\R^n,\;\bm{s}\in\R^{n\times 2},\;\bar{\bm{z}}\in\R_+^{n-1},\;\bm{r}\in\R_+^{n-1},
	\end{aligned}
\end{equation}
where $\lambda_1=1/\sigma_1^2$, $\lambda_n=0$, $\lambda_i=\frac{1}{2\sigma_i^2}$ and $L_i=\lambda_i(1/\sigma_{i+1}^2-\lambda_{i+1})(\tau_{i+1}-\tau_i)+\lambda_i+1/\sigma_{i+1}^2-\lambda_{i+1}$ for $1<i<n$.  For comprehensive details on how this formulation was derived, we direct readers to the original paper by \cite{gomez2021outlier}. The only difference from the formulation in \cite{gomez2021outlier} is that the cardinality constraint in literature is replaced by the penalizing term $c\sum_{i}z_i$ in  \eqref{eq:strong mip}.

\noindent$\bullet$ \sm. In \sm, we solve \eqref{eq:outlier-detection} as a binary submodular minimization problem. Note that 
\eqref{eq:outlier-detection} can be cast in the following matrix form
\begin{equation}\label{eq:outlier-detection-matrix}
	\hspace{-0.5em}\begin{aligned}		
		\min_{\bm{x,z,w}}\;&\frac{1}{2}\bm{x}^\top\bm{P}\bm{x}+\frac{1}{2}(\bm{y}+\bm{w}-\bm{\mu}-\bm{x})^\top\bm{D}(\bm{y}+\bm{w}-\bm{\mu}-\bm{x})+c\sum_{i=1}^nz_i\\
		\text{s.t. }\;&\bm{w}\circ(\ones-\bm{z})=0,\;\bm{z}\in\{0,1\}^n,
	\end{aligned}
\end{equation}
where $\bm{D}$ is a diagonal matrix defined by $D_{ii}=\sigma_i^2$ for $i\in[n]$, and $\bm{P}$ is a Stieltjes matrix given by 
\[ P_{ij}=\begin{cases}
	0&\text{if }j>i+1\\
	-\frac{1}{\tau_{i+1}-\tau_i}&\text{if }j=i+1\\
	\frac{1}{\tau_1}+\frac{1}{\tau_2-\tau_1}&\text{if }i=j=1\\
	\frac{1}{\tau_i-\tau_{i-1}}+\frac{1}{\tau_{i+1}-\tau_i}&\text{if }1<i=j<n\\
	\frac{1}{\tau_n-\tau_{n-1}}&\text{if }i=j=n\\
	P_{ji}&\text{if }i>j.
\end{cases}\]
By minimizing over free variables $\bm{x}$, \eqref{eq:outlier-detection-matrix} can be simplified to
\begin{equation}\label{eq:outlier-detection-sm}
	\begin{aligned}
		\min_{\bm{w,z}}\;&\frac{1}{2}(\bm{w}-\bm{\mu}+\bm{y})^\top\bm{Q}(\bm{w}-\bm{\mu}+\bm{y})+c\sum_{i=1}^nz_i\\
		\text{s.t. }\;&\bm{w}\circ(\ones-\bm{z})=0,\;\bm{z}\in\{0,1\}^n,
	\end{aligned}
\end{equation} 
where $\bm{Q}=\bm{D}-\bm{D}(\bm{P}+\bm{D})^{-1}\bm{D}$ remains a Stieltjes matrix because the inverse of the Stieltjes matrix $P+D$ is componentwise nonnegative. Since \eqref{eq:outlier-detection-sm} is in the form of \eqref{eq:miqp}, it can be solved as a binary submodular minimization problem  which we call \sm. Additionally,  we utilize Algorithm~\ref{algo:active-set} to calculate the extreme bases incurred in the implementation of \sm.  

\subsubsection{Results} 
The dataset used in this study is sourced from the Conic Benchmark Library (CBLIB)\footnote{\url{https://cblib.zib.de/}} \cite{friberg2016cblib}, containing five instances of $(\bm{\tau,\bm\mu,\bm{y},\bm{\sigma}})$ for each $n\in\{100,200,500\}$. The method \sm method is executed on the laptop detailed in Section~\ref{sec:mrf-nonnegative}.  However, to comprehensively evaluate and appreciate the efficiency of the proposed method \sm, the MIP formulations \bigM and \strongMip are executed on high-performance NEOS servers\footnote{\url{https://neos-server.org/neos/}} and solved using Gurobi~10.0.0. Indeed, to solve these formulations, we directly use AMPL files provided by the author of \cite{gomez2021outlier}. A time limit of 1800 seconds is enforced for all three methods. With above setting, the computational results with varying anomaly weight $c$ are shown in Table~\ref{tab:free-MRF}, where the columns Time, Gap, EB and \sparsity are akin in definition those in Section~\ref{sec:mrf-nonnegative-results}. Note that here, \sparsity should be interpreted as the portion of outliers for the robust MRF problem. Each row of the table encapsulates the average performance over five instances. It is worth noting that since not all instances can be solved to optimality within the time limit, \sparsity is solely computed and averaged for the ones solved. 

As one can observe in Table~\ref{tab:free-MRF}, \bigM can solve only five instances to optimality, showcasing the least favorable performance. On the other hand, \strongMip performs better than \bigM -- it is capable of solving 22 instances and achieves notably smaller optimality gaps for those unsolved instances, which is consistent with the results in \cite{gomez2021outlier}. In comparison, \sm is able to solve 90\% of the total 60 instances in a solution time ten times faster than the alternatives\emph{, despite running on a laptop instead of the NEOS server}. Furthermore, we note that besides dimension $n$, \sparsity is another critical factor influencing the performance of both the MIP approach and the submodular minimization approach. As \sparsity increases,  problem \eqref{eq:outlier-detection} becomes more challenging to solve. For instance, \sm can solve all instances with a \sparsity of less than 50 in one minute. Nonetheless, it has difficulty in solving instances when $\sparsity\ge 50$ and $n\ge200$. These challenging scenarios also correspond to a noteworthy increase in solution time. In summary, we ascertain that \sm surpasses existing state-of-the-art MIP approaches, rendering it a favorable choice for addressing \eqref{eq:outlier-detection}.
\begin{table}[h]
	\caption{Results for outlier detection}\label{tab:free-MRF}
	\centering
	\resizebox{\linewidth}{!}{
	\begin{tabular}{cc|rrr|rrr|rcr|c}
		\toprule
		\multirow{2}{*}{$n$} &\multirow{2}{*}{$c$}&\multicolumn{3}{c|}{\bigM}&\multicolumn{3}{c|}{\strongMip}&\multicolumn{3}{c|}{\sm}&\multirow{2}{*}{\sparsity} \\\cmidrule(lr){3-5} \cmidrule(lr){6-8} \cmidrule(lr){9-11}
		& &Time(s) & Gap(\%) & \# & Time(s) & Gap(\%) & \# &   Time(s) &  EB(\%) & \#   & \\\midrule
		
		\multirow{5}{*}{$100$} &        0.1 &     904 &    12.26 &   3 &1800 &11.19 &0 &    148.54 &    41.07 &          5 &         63 \\
		
		&        0.2 &       1800 &     32.05 &0 & 1729& 7.83& 1&      76.64 &    49.68 &          5 &         55 \\
		
		&        0.5 &       1800 &    20.78 &         0 & 1083& 6.24&2&    16.51 &    65.88 &          5 &         36 \\
		
		&          1 &       1470 &    19.25 &          1 & 1081& 3.34&2&    5.56 &    70.03 &          5 &         23 \\
%
		\bottomrule
		\multirow{5}{*}{$200$} &        0.1 &       1800 &    57.51 &          0 & 1800 & 24.07 &0 &     904.96 &    58.49 &          3 &         62 \\
		
		&        0.2 &       1800 &    55.19 &          0 & 1483& 15.70& 1&   544.50 &    72.47 &          4 &         51 \\
		
		&        0.5 &       1800 &   47.68 &          0 & 773& 6.31&3&  54.37 &    89.79 &          5 &         32 \\
		
		&          1 &       1475 &    31.27 &          1 & 721& 1.85&3&  11.36 &    93.98 &          5 &         14 \\
		
		\bottomrule
		\multirow{5}{*}{$500$} &        0.1 &       1800 &    77.38 &          0 & 1800&18.50&0&    1193.66 &    95.19 &          2 &         50 \\
		
		&        0.2 &       1800 &    70.89 &          0 & 1523& 8.36&1&  539.86 &    98.64 &          5 &         50 \\
		
		&        0.5 &       1800 &    53.59 &          0 & 710&0.54&4&   42.86 &    99.81 &          5 &         11 \\
		
		&          1 &       1800 &    29.83 &          0 & 6& 0.00&5&    12.93 &    99.85 &          5 &          3 \\
		
		\bottomrule
	\end{tabular}  
}
\end{table}

\section{Conclusion}\label{sec:conclusion}
In this paper, we study a class of convex submodular minimization problems with indicator variables, of which the inference of Markov random fields with sparsity and robustness priors is a special case. Such a problem can be solved as a binary submodular minimization problem and thus in (strongly) polynomial time provided that for each fixed binary variable, the resulting convex optimization subproblem is (strongly) polynomially solvable. When applied to quadratic and conic quadratic cases, it extends known results in the literature.  More efficient implementations are also proposed by exploiting the isotonicity of the solution mapping in parametric settings.

\section*{Acknowledgments}
The authors thank Professor Jong-Shi Pang for his valuable discussion and suggestions during the development of this work.

Andr\'es G\'omez is supported, in part, by the Air Force Office of Scientific Research under grant No. FA9550-24-1-0086. Shaoning Han is supported by the Ministry of Education, Singapore, under the Academic Research Fund Tier 1 (FY2024).

\bibliographystyle{apalike}
\bibliography{Bibliography}

\appendix
\section{Cutting plane method for binary submodular function minimization}\label{sec:cutting-plane}
Given a binary submodular function $v:\set{Z}\to\R$, where $\set{Z}\subseteq\{0,1\}^n$ is a lattice, we aim to solve $\min\limits_{\bm{z}\in\set{Z}}\,v(\bm z)$. Without loss  of generality we assume that $v(\bm{0})=0$; otherwise, one can consider the function $v(\bm{z})-v(\bm{0})$. For any vector $\bm{\bar z}\in[0,1]^n$, define function $v_L(\bm z;\bm{\bar z})\defeq \sum_{i=1}^n\left( v\left(\ones^{\pi_{[i]}}\right)-v\left(\ones^{\pi_{[i-1]}}\right)\right)z_{\pi_i}$, where $\pi\in\Pi([n])$ such that $\bar{z}_{\pi_1}\ge\bar{z}_{\pi_2}\dots\ge\bar{z}_{\pi_n}$.   Note that the Lov\'{a}sz extension of $v(\cdot)$ can be expressed as $\displaystyle v^L(\bm{z})=\max_{\bm{\bar z}\in[0,1]^n}v_L(\bm{z};\bm{\bar z})$  which is actually the maximum of a finite (but exponential in $n$) number of linear functions  \cite{lovasz1983submodular}. Moreover, ${v}^L(\bm{\bar z})=v_L(\bm{\bar z};\bm{\bar z})$ holds for all $\bm{\bar z}\in[0,1]^n$. Since $\displaystyle \min_{\bm{z}\in\{0,1\}^n}v(\bm{z})=\min_{\bm{z}\in[0,1]^n}{v}^L(\bm{z})$ is equivalent to a linear program with an exponential number of constraints
\begin{align*}
	 \min_{(t,\bm{z})\in\R^{n+1}}\;\;&t\\
	 \text{s.t. }\;& t\ge v_L(\bm{z};\bm{\bar z})\quad\forall \bm{\bar z}\in[0,1]^n,
\end{align*} the submodular function minimization problem can be solved using the standard cutting plane method, where according to the touching property ${v}^L(\bm{\bar z})=v_L(\bm{\bar z};\bm{\bar z})$, the separating oracle is induced by sorting the elements of the incumbent solution $\bm{\bar z}$; see \cite{atamturk2022submodular} or Section~6.3 in \cite{bertsimas1997introduction} for details.

\end{document}